\documentclass[12pt,draftcls,onecolumn]{IEEEtran}

\ifCLASSINFOpdf
\else
\fi

\usepackage{graphics} 
\usepackage{times} 
\usepackage{amsmath} 
\usepackage{amssymb}  
\usepackage{amsfonts}
\usepackage{subfig}
\usepackage{cite}
\usepackage[colorlinks=true]{hyperref}
\usepackage{bm}
\usepackage{multirow}
\usepackage{comment}
\usepackage{mathtools}
\usepackage{tikz}
\usepackage{url,bm,xspace,dsfont}
\usepackage{cases}                                        
\usepackage{ntheorem}
\usepackage{verbatim}
\usepackage{proof}
\usepackage{tikz}
\usepackage{soul}

\usetikzlibrary{matrix,arrows,positioning}
\newtheorem{theorem}{Theorem}[section]
\newtheorem{lemma}[theorem]{Lemma}
\newtheorem{definition}[theorem]{Definition}

\newtheorem{remark}[theorem]{Remark}

\newtheorem{example}[theorem]{Example}

\newtheorem{problem}[theorem]{Problem}
\newtheorem{algorithm}[theorem]{Algorithm}

\newcommand{\sr}{\stackrel}

\newcommand{\tri}{\sr{\triangle}{=}}

\newcommand{\be}{\begin{equation}}
\newcommand{\ee}{\end{equation}}
\newcommand{\bea}{\begin{eqnarray}}
\newcommand{\eea}{\end{eqnarray}}
\newcommand{\bes}{\begin{eqnarray*}}
\newcommand{\ees}{\end{eqnarray*}}
\newcommand{\bi}{\begin{itemize}}
\newcommand{\ei}{\end{itemize}}
\newcommand{\ben}{\begin{enumerate}}
\newcommand{\een}{\end{enumerate}}
\newcommand{\bp}{\begin{problem}}
\newcommand{\ep}{\end{problem}}
\newcommand{\hso}{\hspace{.1in}}
\newcommand{\hst}{\hspace{.2in}}

\newcommand{\noi}{\noindent}

\def\supp{\mathop{\mathrm{supp}}}  

\newcommand{\calX}{\mathcal{X}}
\newcommand{\calY}{\mathcal{Y}}

\begin{document}
%
\title{\bf Approximation of Markov Processes by Lower Dimensional Processes via\\ Total Variation Metrics}

\author{Ioannis~Tzortzis, Charalambos~D.~Charalambous, Themistoklis~Charalambous, \\ Christoforos N. Hadjicostis and Mikael Johansson
\thanks{I. Tzortzis, C. D. Charalambous and C. N. Hadjicostis are with the Department of Electrical Engineering, University of Cyprus, Nicosia, Cyprus.
      Emails:  {\tt \{\small tzortzis.ioannis, chadcha,chadjic\}@ucy.ac.cy}.}%
\thanks{T. Charalambous and M. Johansson are with the School of Electrical Engineering, Royal Institute of Technology (KTH), Stockholm, Sweden.
      Email:  {\tt\small \{themisc, mikaelj\}@kth.se}.}%
      }
\maketitle
%
%
%
%
%
\begin{abstract}
The aim of this paper is to approximate a finite-state Markov process by another process with fewer states, called herein the approximating process. The approximation problem is formulated using two different methods. 

The first method, utilizes the total variation distance to discriminate the transition probabilities of a high dimensional Markov process and a reduced order Markov process. The approximation is obtained by optimizing a linear functional defined in terms of transition probabilities of the reduced order Markov process over a total variation distance constraint. The transition probabilities of the approximated Markov process are given by a water-filling solution.

The second method, utilizes total variation distance to  discriminate the invariant probability of a Markov process and that of the approximating process. The approximation is obtained via two alternative formulations: (a) maximizing a  functional of the occupancy distribution of the Markov process, and (b) maximizing the entropy of the approximating process invariant probability. For both formulations, once the reduced invariant probability is obtained, which does not correspond to a Markov process, a further approximation by a Markov process is proposed which minimizes the Kullback-Leibler divergence. These approximations are given by water-filling solutions.

Finally, the theoretical results of both methods are applied to specific examples to illustrate the methodology, and the water-filling behavior of the approximations.
\end{abstract}


\IEEEpeerreviewmaketitle

%
%
%
%
\section{Introduction}\label{sec.Intro}

Finite-State Markov (FSM) processes are often employed to model physical phenomena in many diverse areas, such as machine learning, information theory (lossy compression), networked control and telecommunication systems, speech processing, systems biology, etc. In many of these  applications the state-space of the Markov process is prohibitively large, to be used in analysis and  simulations. One approach often pursue to overcome the large number of states is to  approximate the Markov process by a lower dimensional Markov process, with respect to certain measures of discriminating or approximating the distribution of the high dimensional Markov process by a reduced one. Such methods are described using relative entropy as a measure of approximation in  \cite{Vidyasagar:10,Beck:11,Deng:11,2012:SharmaTAC} (and references therein). Further discussion of model reduction methods for Markov chains can be found in \cite{beck:2009}. In general, approximating a Markov process by another process subject to a fidelity of reproduction is not  necessarily Markov, but a finite-state hidden Markov process. This is a well known result of Information Theory \cite{2006:Cover}, on lossy compression of Markov sources with respect to a fidelity criterion. Model reduction of hidden Markov models via aggregation can be found in \cite{deng:2010,Vidyasagar:10,DBLP:conf/cdc/DengMMV11}. Specifically, in \cite{DBLP:conf/cdc/DengMMV11} the aggregated hidden Markov model is expressed as a function of a partition function and a recursive learning algorithm is proposed, which solves the optimal partition problem.

In this paper, the approximation problem of a FSM process by another process (FSM or FSHM) with reduced state-space is formulated as an optimization problem, with respect to a certain pay-off subject to a fidelity criterion defined by the total variation distance metric, using two different methods which are elaborated below.\\


\noi\textit{Method 1. }
 
Approximate the transition probabilities of a FSM process by another FSM process with reduced transition probability matrix. This approximation problem is formulated as a maximization of a linear functional on the transition probabilities of the reduced FSM processes, subject to a fidelity criterion defined by the total variation distance between the transition probabilities of the high and low FSM process. The main contributions of this method are the following:
\ben
\item[(i)] a direct method for Markov by Markov approximation based on the transition probabilities of the original FSM process, exhibiting a water-filling behavior;
\item [(ii)] an example which illustrates the methodology, and the properties of the approximation.\\
\een

\noi\textit{Method 2.}

Approximate a FSM process by another process with lower dimensional state-space, without imposing the assumption that the approximating process is also a Markov process. The following two formulations are investigated:
\begin{itemize}
\item[(a)] maximize an average pay-off, described in terms of the occupation measure of the high dimensional Markov process, subject to a fidelity criterion defined by the total variation distance metric, between the invariant distribution of the higher dimensional Markov process and that of  the lower dimensional process.
 
\item[(b)] maximize the entropy (Jayne's maximum entropy \cite{Jaynes:1957}) of the invariant distribution of the lower dimensional process, subject to a fidelity criterion defined by the total variation distance metric, between the invariant distribution of the higher dimensional Markov process and that of the  the lower dimensional process. 
\end{itemize}
\noi For both formulations, the resulting approximated process is not necessarily Markov, but a hidden Markov process. The crux of the approach considered lies in finding an optimal partition function which aggregates states of the original FSM process to form the reduced order process. Moreover, an approach is described to further approximate the hidden Markov process by a Markov process, by minimizing the Kullback-Leibler divergence. The main contributions of this method are the following:
\ben
\item[(i)] iterative algorithms to compute the invariant distribution of the approximating process;
\item[(ii)] extremum measures which exhibit water-filling behavior, and solve the approximation problems;
\item[(iii)] optimal partition functions which aggregate the original FSM process to form the reduced order processes;
\item[(iv)] examples which illustrate the approximation method and the properties of solutions to both formulations.
\een

The rest of the paper is organized as follows. In Section~\ref{sec.Problem.Form}, the total variation distance and the Kullback-Leibler divergence rate are defined, and the approximation problems are introduced. In Section~\ref{sec.sol.opt.problem.2}, the solution of approximation problem based on \emph{Method 1} is given. In Sections~\ref{subsec.occ.problem.} and \ref{subsec.entr.problem.}, the solution of approximation problems based on \emph{Method 2} is given.  In Section~\ref{sec.examples}, several examples are presented to illustrate the approximation methods. Section~\ref{sec.conclusions} concludes by discussing the most important results obtained in this paper.

%
%
%
%
\section{Problem Formulation}\label{sec.Problem.Form}

%
%
\subsection{Preliminaries and discrepancy measures}\label{subsec.preliminaries}
We consider a discrete-time homogeneous Markov process $\{X_t:t=0,1,\dots\}$, with state-space $\calX$ of finite cardinality  $card(\calX)=|\calX|$, and transition probability matrix $P$ with elements $\{p_{ij}:i,j=1,\dots,|\calX|\}$ defined by
\begin{equation*}
p_{ij}\tri \mathbb{P}(X_{t+1}=j|X_t=i),\quad i,j\in\calX,\quad t=0,1,\dots .
\end{equation*}
The Markov process is assumed to be irreducible, aperiodic having a unique invariant distribution $\mu=[\mu_1\ \mu_2\dots\mu_{|\calX|}]$ satisfying 
\begin{equation*}
\mu=\mu P .
\end{equation*}
For the rest of the paper we adopt the notation $(\mu,P,\calX)$ to denote a stationary FSM process, with transition probabity matrix $P$, stationary distribution $\mu$, and state-space $\calX$.

The distance metrics we will use to define the discrepancy between two probability distributions (and conditional probability distributions) are the Total Variation distance, and the Kullback-Leibler divergence. These are introduced below.

Consider the finite alphabet space $(\calX,\cal M)$, with ${\cal M}={2}^{|\calX|}$. Define the set of probability vectors on $\calX$ by 
\begin{equation*}
\mathbb{P}(\calX)\tri \Big\{p=(p_1,\dots,p_{|\calX|}):p_i\geq 0,i\in\calX,\sum_{i\in\calX}p_i=1\Big\} .
\end{equation*}
Thus, $p\in\mathbb{P}(\calX)$ is a probability vector in $\mathbb{R}_+^{|\calX|}$.

\subsubsection{Total Variation (TV) distance}\label{subsubsec.TV} \cite{dunford} The TV distance is a metric $||\cdot||_{TV}:\mathbb{P}(\calX)\times \mathbb{P}(\calX)\longrightarrow [0,2]$ defined by
\begin{equation*}
||\nu-\mu||_{TV}\tri \sum_{i\in\calX}|\nu_i-\mu_i|,\qquad \nu,\mu\in \mathbb{P}(\calX).
\end{equation*}

\subsubsection{Relative Entropy distance}\label{subsubsec.RE}
\cite{Dupuis97} The relative entropy  of $\nu\in \mathbb{P}(\calX)$ with respect to $\mu\in \mathbb{P}(\calX)$ is a mapping ${\mathbb D}({ \cdot}|| { \cdot} ):\mathbb{P}(\calX)\times \mathbb{P}(\calX)\longrightarrow [0,\infty]$ defined by
\bes
{\mathbb D}({ \nu}|| { \mu} ) \tri \sum_{i\in\calX}\nu_i\log\frac{\nu_i}{\mu_i}.
\ees
 It is well known that ${\mathbb D}({ \nu}|| { \mu} )\geq0,\forall\nu,\mu\in \mathbb{P}_1(\calX)$, while ${\mathbb D}({ \nu}|| { \mu} )=0\Leftrightarrow\nu=\mu$.

Given a probability vector $\mu\in \mathbb{P}(\calX)$ define the fidelity set via the ball, with respect to the TV distance, centered at the vector $\mu\in \mathbb{P}(\calX)$, having radius $R\in[0,2]$ by 
\begin{equation} \label{TVClass}
{\mathbb B }_R({ \mu})   \tri \Big\{  \nu \in \mathbb{P}(\calX): ||{ \nu} -{ \mu}||_{TV} \leq R \Big\}.
\end{equation}
The two extreme cases of this set are $R=0$ implying $\nu_i=\mu_i$, $\forall i\in\calX$, $a.e.$, and $R=2$ implying that the support sets of ${ \nu}$ and ${ \mu}$ denoted by $\supp({ \nu})$ and   $\supp({ \mu})$, respectively, are non-overlapping, that is, $\supp({ \nu}) \cap \supp({ \mu}) = \emptyset$. One of the most interesting properties of TV distance ball is that, any probability vector $\nu\in {\mathbb B }_R({ \mu}) $ may not be absolutely continuous with respect to ${ \mu}$ (i.e., $\mu_i =0$ for some $i \in {\calX}$ then $\nu_i=0$). Consequently, any approximating probability vector $\nu\in {\mathbb B }_R({ \mu})$ can be defined on an alphabet $\calY$ with smaller cardinality than the probability vector ${ \mu}\in \mathbb{P}(\calX)$, that is, $\supp(\nu)\subseteq \supp(\mu)$. The total variation metric is also discussed in \cite{Vidyasagar12}.
There is an anthology of distances and distance metrics on the space of probability distributions which are related to total variation distance \cite{gibbs}, and therefore one can obtain various lower and upper bounds on the performance with respect to other types of discrepancy measures. For example, if $ { \nu}$ is not absolutely continuous with respect to ${ \mu}$, $\forall \nu\in {\mathbb B }_R({ \mu})$ belonging to total variation distance class, by Pinsker's inequality \cite{1964:Pinsker}, then
\bes
 ||{ \nu} - { \mu}||_{TV}^2 \leq
2 {\mathbb D}( { \nu} || { \mu} ), \quad \forall \nu\in {\mathbb B }_R({ \mu}),\ { \mu} \in \mathbb{P}_1(\calX).
 \ees
This is one such relation between $||\cdot||_{TV}$ and ${\mathbb D}(\cdot|| \cdot)$.

Let $(\mu,P,\calX)$ and $(\nu,\Phi,\calX)$ be two stationary FSM processes. A version of the KL divergence used in \cite{2006:Cover,Rached:04}, is defined by
\begin{equation}\label{eq.KL}
\mathbb{D}_{\mu}(P||\Phi)\tri\sum_{i,j\in \calX}\mu_iP_{ij}\log\Big(\frac{P_{ij}}{\Phi_{ij}}\Big) ,
\end{equation}
where $P_{i\bullet}$ is assumed to be absolutely continuous with respect to $\Phi_{i\bullet}$, that is, for any $i\in\calX$, $\Phi_{ij}=0$ for some $j\in\calX$ then $P_{ij}=0$. Note that \eqref{eq.KL} is used to compare stationary Markov processes which are defined on the same state-space. For Markov processes which are defined on different state-spaces, \eqref{eq.KL} is defined with respect to the lifted version of the lower dimensional Markov process (see \cite{Deng:11}), defined by
\begin{equation}\label{eq.lifted.chain}
\widehat{\Phi}_{ij}=\frac{\mu_j}{\sum_{k\in\psi(j)}\mu_k}\Phi_{\varphi(i)\varphi(j)},\qquad i,j\in\calX ,
\end{equation}
where $\psi(j)$ denotes the set of states belonging to the same group as the $j$th state, and $\varphi$ denotes a partition function from $\calX$ onto $\calY$. For the rest of the paper we will use the notation $\mathbb{D}^{(\varphi)}(P||\Phi)=\mathbb{D}_\mu(P||\widehat{\Phi})$ to denote the KL divergence distance between two Markov processes via liflting.

%
%
\subsection{Approximation problems}\label{subsec.approx.problems}
In this section we introduce the approximation problems described in the introduction. We propose two different methods to approximate FSM processes by lower dimensional processes, as follows.\\

\subsubsection{Method 1}\label{subsubsec.approx.method2}
This method is based on comparing two FSM processes $(\mu,P,\calX)$ and $(\nu,\Phi,\calY)$, $\calY\subseteq \calX$, by working directly on their transition probability matrices $P$ and $\Phi$. The approximation problem is formulated as a maximization of a linear functional, defined on the transition probabilities of the reduced order FSM process $(\nu,\Phi,\calY)$, subject to a TV  distance fidelity criterion, between the transition probabilities of the high and low dimensional FSM processes. The precise problem formulation is given below.
\begin{problem}\label{problem3}
Given a FSM process $(\mu,P,\calX)$, find a transition probability matrix $\Phi$ which solves the maximization problem defined by
\begin{align}\label{MCproblem1}
\max_{\Phi_{i\bullet}\in\mathbb{P}(\calY),\forall i\in\calY} & \quad\sum_{i\in \calX}\sum_{j\in \calX}\ell_j \Phi_{ij}\mu_i\\
\mbox{s.t.} & \quad\sum_{i\in\calX}\sum_{j\in\calX}|\Phi_{ij}-P_{ij}|\mu_i\leq R,\quad \forall R\in[0,2]\nonumber.
\end{align}
where $\ell\tri \{\ell_1,\dots,\ell_{|\calX|}\}\in\mathbb{R}_+^{|\calX|}$ (i.e., set of non-negative vectors of dimension $|\calX|$).
\end{problem}
The choice of $\ell$ weights the transition probabilities. 

The optimal transition probability matrix $\Phi$ which solves maximization problem \eqref{MCproblem1} is obtained for all values of TV parameter $R\in[0,2]$, and exhibits a \textit{water-filling} solution. In addition, as the TV parameter increases, it turns out that the dimension of the transition matrix $\Phi$  is reduced, and hence, a reduced order FSM process is obtained.\\

\subsubsection{Method 2}\label{subsubsec.approx.method1}
Given a FSM process $(\mu,P,\calX)$ and a parameter $R\in[0,2]$, define the average pay-off with respect to the stationary distribution $\nu\in \mathbb{B}_R(\mu)\subset \mathbb{P}(\calX)$ by 
\begin{equation}\label{lin.funct.}
\mathbb{L}(\nu)=\sum_{i\in\calX}\ell_i\nu_i, \quad \ell\in\mathbb{R}_+^{|\calX|}.
\end{equation}
The objective is to approximate $\mu\in\mathbb{P}(\calX)$ by $\nu\in\mathbb{B}_R(\mu)$, by solving the maximization problem defined by 
\begin{equation}\label{max.aver.equation}
\mathbb{L}(\nu^*)=\max_{\sr{\nu\in\mathbb{B}_R(\mu)}{\mu=\mu P}}\mathbb{L}(\nu),\quad \forall R\in[0,2], 
\end{equation}
for two alternative choices of the parameters $\ell\in \mathbb{R}_+^{|\calX|}$, as follows.

 \textbf{Formulation (a)} (\textit{Approximation Based on Occupancy Distribution})

\noi Let $\ell_i\triangleq \mu_i$, $\forall i\in\calX$, which implies \eqref{max.aver.equation} is equivalent to maximizing a weighted sum of the stationary distribution $\{\nu_i:i\in\calX\}\in \mathbb{P}(\calX)$, subject to a fidelity criterion. This formulation leads to an approximation algorithm described via reduction of the states (i.e., by deleting certain states of the original Markov process) to  obtain the approximating reduced state process. Intuitively, the optimal solution has the property of maintaining and strengthening the states with the highest invariant probability, while removing the states with the smallest invariant probability. 

\textbf{Formulation (b)} (\textit{Approximation Based on Maximum Entropy Principle})

\noi Let $\ell_i\triangleq -\log\nu_i$, $\forall i\in\calX$, which implies that \eqref{max.aver.equation} is equivalent to the problem of finding the approximating distribution corresponding to the minimum description codeword length \cite{Barron:1998}. This formulation leads to an optimal approximation algorithm described via aggregation of the states (i.e., by grouping certain states of the original Markov process)  to  obtain the approximated reduced state process, which is a hidden Markov process. This formulation  is related to minimizing the average codeword length of the approximated Markov process, subject to a fidelity criterion.\\

The approximated probability vector is based on the following concept. Given a FSM process $(\mu,P,\calX)$, the optimal probabilities of the reduced process are defined on $\calX$, which is partitioned into disjoint sets $\calX=\cup_{i=1}^K\calX_i$, $K\leq |\calX|$. The solution of the optimization problems based on Method 2(a) and 2(b), give the maximizing probability $\nu^*(\calX_i)$, $i=1,\dots,K$, on this partition. 

For Method $2(a)$, as $R$ increases the maximizing probability vector, $\nu^*$, is given by a \textit{water-filling} solution, having the property that states of the initial probability vector $\mu\in \mathbb{P}(\calX)$ are deleted to form a new partition of $\calX$, denoted by $\calX=\cup_{i=1}^M\calY_i$, $M\leq K\leq |\calX|$. The approximated probability vector is then obtained as defined below.
\begin{definition}\label{lower.dim.process1}(Approximated Probability Vector based on Occupancy Distribution) \\Define the restriction of $\nu^*$ on only those elements of the partition $\{\calY_1,\dots,\calY_M\}$ which have non-zero probability by
\begin{equation} \nu^*|_{\supp(\nu^*)\neq 0}:\{\calY_{i_1},\calY_{i_2},\dots,\calY_{i_k}\}\longmapsto [0,1],\end{equation}  
where $\{\calY_{i_1},\calY_{i_2},\dots,\calY_{i_k}\}\subseteq \{\calY_1,\dots,\calY_M\}$, and $i_1,i_2,\dots,i_k\in \{1,2,\dots,M\}$. The approximated probability vector based on occupancy distribution is defined by 
\begin{equation}
\bar{\nu}=\nu^*|_{\supp(\nu^*)\neq 0} ,
\end{equation} 
having states which are in one-to-one correspondence with $\{1,2,\dots,k\}$, via the mapping $\calY_{i_1}\longmapsto 1$, $\calY_{i_2}\longmapsto 2,\dots$, $\calY_{i_k}\longmapsto k$, with corresponding process $\{Y_t:t=0,1,\dots\}$ having state-space $\calY=\{1,2,\dots,k\}$.
\end{definition}

For Method $2(b)$, as $R$ increases the maximizing probability vector $\nu^*$, exhibits a \textit{water-filling} solution, with the property that states of $\mu\in\mathbb{P}(\calX)$ are aggregated together to form a new partition of $\calX$. The approximated probability vector is obtained as defined below.
\begin{definition}\label{lower.dim.process2}(Approximated Probability Vector based on Maximum Entropy Principle)\\
Define $\bar{\nu}=\nu^*$ if all elements of $\nu^*(\calX_k)$ are not equal and the state-space of $\bar{\nu}$ is $\calY=\{1,\dots,K\}$. If any of the $\nu^*(\calX_k)$, $k\in\{1,\dots,K\}$ become equal then a new probability vector $\bar{\nu}$ is defined by adding together those $\nu^*\in\mathbb{P}(\calX)$ which are equal, and setting $\bar{\nu}\tri\nu^*(\calX_k)$ for the $\nu^*(\calX_k)$ whose elements are not equal. The resulting approximated probability vector based on maximum entropy principle $\bar{\nu}\in \mathbb{P}(\calY)$, with corresponding process $\{Y_t:t=0,1,\dots\}$, is defined on a state-space $\calY$, whose cardinality is less or equal to $|\calX|$. 
\end{definition}

\begin{remark}
In general, the reduction based on Methods $2(a),(b)$ do not lead to a Markov chain, even though it could be the case.
\end{remark}

However, an Markov approximating process is obtained by the following two-step procedure. Step 1 corresponds to the  the approximating problems described above\footnote{The reduced approximating process is obtained without a priori imposing the assumption that it is also a Markov process.}. Step 2 utilizes the approximating process $\{Y_t:t=0,1,\dots\}$ of step 1,  to further approximate a FSM process by another FSM process $(\bar{\nu},\Phi,\calY)$, $\calY\subset\calX$. Here, the objective is to find an optimal partition function $\varphi$ and a transition matrix $\Phi$ which minimizes the KL divergence rate \cite{Deng:11} defined by
\begin{equation}\label{eq.KL.lifted}
\mathbb{D}^{(\varphi)}(P||\Phi)=\sum_{i,j\in \calX}\mu_iP_{ij}\log\Big(\frac{P_{ij}}{\widehat{\Phi}_{ij}}\Big) ,
\end{equation}
where $\widehat{\Phi}$ is given by \eqref{eq.lifted.chain}, and denotes the lifted version of the lower dimensional Markov chain $\Phi$ by using an optimal partition function $\varphi$. By employing certain results from \cite{Deng:11}, the transition matrix $\Phi$ which solves \eqref{eq.KL.lifted} is obtained. What remains, is to find an optimal partition function $\varphi$, for the approximation problems of Method $2(a)$ and $2(b)$. This Markov by Markov approximation is found by working only with values of TV parameter for which a reduction of the states occurs, that is, $|\calY|<|\calX|$.

\begin{figure}[!h]
\centering
\includegraphics[width=1\columnwidth]{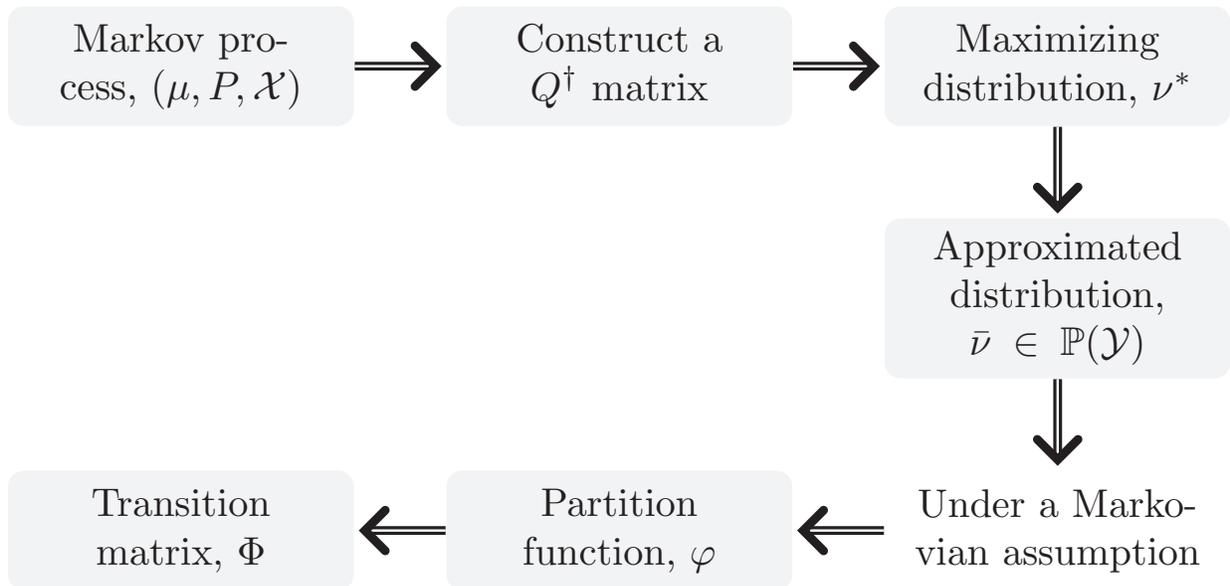}
\caption[]{Procedure of \emph{Method 2}.}\label{fig.proc}
\end{figure}

Given a FSM process $(\mu,P,\calX)$, an algorithm is presented, which describes how to construct the transition probability matrix $Q^{\dagger}$, from the maximizing distribution $\nu^*$ of problem \eqref{max.aver.equation} for Method $2(a)$ and $2(b)$. Then, using Definitions \ref{lower.dim.process1} and \ref{lower.dim.process2}, a lower probability distribution $\bar{\nu}\in \mathbb{P}(\calY)$ is obtained. Under the restriction that the lower dimensional process is also a FSM process $(\bar{\nu},\Phi,\calY)$, $\calY\subset \calX$, an optimal partition function $\varphi$ and a transition probability matrix $\Phi$, are found which minimize the KL divergence rate between $P$ and $\widehat{\Phi}$. The approximation procedure for Method $2(a)$ and $2(b)$, is shown in Fig.\ref{fig.proc}.

The precise problem definition of approximation Method 2 based on occupancy distribution is given below.

\begin{problem}\label{problem.occ}(Approximation Based on Occupancy Distribution)\\ Let $\{\ell_i:i\in\calX\}\in \mathbb{R}_+^{|\calX|}$ denote the occupancy distribution of a FSM process $(\mu,P,\calX)$ defined by $\ell_i\triangleq \mu_i$, $\forall i\in\calX$. Find $\{\nu_i:i\in\calX\}\in \mathbb{P}(\calX)$ which solves
\begin{equation}\label{max.occupancy}
\max_{\sr{\nu\in {\mathbb B }_R({ \mu})}{\mu=\mu P}}\sum_{i\in\calX}\mu_i\nu_i.
\end{equation} 
Given the optimal solution of \eqref{max.occupancy}, let $\bar{\nu}$ of Definition \ref{lower.dim.process1} denote the invariant distribution of a lower dimensional FSM process $(\bar{\nu},\Phi,\calY)$, $\calY\subset \calX$. 

\noi Find an optimal partition function $\varphi$, and calculate the transition probability matrix $\Phi$, which satisfies $\bar{\nu}=\bar{\nu}\Phi$, and minimizes the KL divergence rate defined by 
\begin{equation}\label{KLoccupancy}
\min_{\sr{\varphi,\Phi}{\bar{\nu}=\bar{\nu}\Phi}}\mathbb{D}^{(\varphi)}(P||\Phi) .
\end{equation}
\end{problem}

Other reasonable choices, are possible by letting $\ell\in\mathbb{R}_+^{|\calX|}$ correspond to a reward or a profit, a cost or a loss, etc., whenever a node is visited.

 Next, the precise problem definition of approximation Method 2 based on maximum entropy principle is given.
\begin{problem}\label{problem.entr}(Approximation Based on Maximum Entropy Principle)\\
Maximize the entropy of $\{\nu_i:i\in\calX \}\in \mathbb{P}(\calX)$ subject to total variation fidelity set, defined by 
\begin{equation}\label{max.entropy}
\max_{\sr{\nu\in {\mathbb B }_R({ \mu})}{\mu= \mu P}}H(\nu), \hst H(\nu) \tri -\sum_{i \in \calX} \log (\nu_i) \nu_i .
\end{equation}
Given the optimal solution of \eqref{max.entropy}, let $\bar{\nu}$ of Definition \ref{lower.dim.process2} denote the invariant distribution of a lower dimensional Markov process $(\bar{\nu},\Phi,\calY)$, $\calY\subset \calX$. 

\noi Find an optimal partition function $\varphi$, and calculate the transition probability matrix $\Phi$, which satisfies $\bar{\nu}=\bar{\nu}\Phi$, and minimizes the KL divergence rate defined by 
\begin{equation}\label{KLentropy}
\min_{\sr{\varphi,\Phi}{\bar{\nu}=\bar{\nu}\Phi}}\mathbb{D}^{(\varphi)}(P||\Phi).
\end{equation}
\end{problem}

Problem \eqref{max.entropy} is of interest when the concept of insufficient reasoning (e.g., Jayne's maximum entropy principle\footnote{The maximum entropy  principle states that, subject to precisely stated prior data, the probability distribution which best represents the current state of knowledge is the one with largest entropy.} \cite{Jaynes:1957}) is applied to construct a model for $\nu\in {\mathbb P}({\calX})$, subject to information quantified via the fidelity set defined by the  variation distance between $\nu$ and  $\mu$.

It is not difficult to show that the maximum entropy approximation defined by \eqref{max.entropy} is precisely equivalent to the problem of finding the approximating distribution corresponding to the minimum description codeword length, also known as the universal coding problem \cite{Barron:1998,Rissanen:1978}, as follows.  Let  $\{\ell_i:i\in\calX\}\in \mathbb{R}_+^{|\calX|}$ denote the positive codeword lengths corresponding to each symbol of the approximated distribution, which satisfy the Kraft inequality of lossless Shannon codes $\sum_{i\in\calX}D^{-\ell_i}\leq 1$, where the codeword alphabet is $D$-ary (unless specified otherwise $\log(\cdot)\tri \log_D(\cdot))$. Then, by the Von-Neumann's theorem, which holds due to compactness and convexity of the constraints, it follows that
\begin{equation*}
\min_{\ell\in \mathbb{R}_+^{|\calX|}:\sum_{i\in\calX}D^{-\ell_i}\leq 1}\max_{\sr{\nu\in {\mathbb B }_R({ \mu})}{\mu=\mu P}}\sum_{i\in\calX}\ell_i\nu_i
=\max_{\sr{\nu\in {\mathbb B }_R({ \mu})}{\mu=\mu P}}\min_{\ell\in \mathbb{R}_+^{|\calX|}:\sum_{i\in\calX}D^{-\ell_i}\leq 1}\sum_{i\in\calX}\ell_i\nu_i
=\max_{\sr{\nu\in {\mathbb B }_R({ \mu})}{\mu=\mu P}} H(\nu).
\end{equation*}
Hence, for $\ell_i\tri -\log\nu_i$, $\forall i \in \calX$, the optimization \eqref{max.aver.equation} is equivalent to optimization \eqref{max.entropy}. 


%
%
%
%
\section{{Method 1:} Solution of approximation problem}\label{sec.sol.opt.problem.2}
In this section, we give the main theorem which characterizes the solution of Problem \ref{problem3}. Define the maximum and minimum values of the sequence $\{\ell_1,\dots,\ell_{|\calX|}\}\in \mathbb{R}_+^{|\calX|}$ by 
\begin{equation*}\ell_{\max}\triangleq \max_{i\in\calX}\ell_i,\qquad \ell_{\min}\triangleq \min_{i\in\calX}\ell_i\end{equation*}
and its corresponding support sets by
\begin{equation*}
\calX^0 \triangleq \{i\in \calX:\ell_i=\ell_{\max} \},\qquad \calX_0\triangleq\{i\in \calX:\ell_i=\ell_{\min} \}.
\end{equation*}
For all remaining elements of the sequence, $\{\ell_i:i\in \calX \setminus \calX^0\cup\calX_0\}$, define recursively the set of indices for which $\ell$ achieves its $(k+1)$th smallest value by $\calX_k$,
where $k\in\{1,2,\hdots, |\calX\setminus \calX^0\cup\calX_0|\}$, till all the elements of $\calX$ are exhausted (i.e., $k$ is at most $|\calX\setminus\calX^0\cup\calX_0|$), and the corresponding values of the sequence on the $\calX_k$ sets by $\ell(\calX_k)$.

 For a fixed $i\in\calX$, define the total variation of a finite signed measure $\Xi_{ij}\tri \Phi_{ij}-P_{ij}$, $\forall j\in \calX$, to be equal to the summation of its positive and its negative part, that is, \begin{equation}\label{tv1}
||\Xi_{i\bullet}||_{TV}\tri \sum_{j\in\calX}\Xi^+_{ij}+\sum_{j\in\calX}\Xi^-_{ij},\qquad \forall i\in\calX .
\end{equation}
By utilizing the fact that $\sum_{j\in\calX}\Xi_{ij}=0$, $\forall i\in\calX$ then \begin{equation}
\sum_{j\in\calX}\Xi^+_{ij}=\sum_{j\in\calX}\Xi^-_{ij}=\frac{||\Xi_{i\bullet}||_{TV}}{2},\qquad \forall i\in\calX.
\end{equation}
Let $\alpha_i\tri ||\Xi_{i\bullet}||_{TV}$, $\forall i\in \calX$, then the constraint of \eqref{MCproblem1} is equivalent to \begin{equation}\label{tvtvconstr}
\sum_{i\in\calX}\alpha_i\mu_i\leq R.
\end{equation} 
and the pay-off can be reformulated as follows.
\begin{equation}\label{eq.payoff2}
\max_{\Phi_{i\bullet}\in \mathbb{P}_1(\cdot)} \sum_{i\in \calX}\sum_{j\in \calX}\ell_j \Phi_{ij}\mu_i\equiv  \sum_{i\in \calX}\sum_{j\in \calX}\ell_j P_{ij}\mu_i+\max_{\Phi_{i\bullet}\in \mathbb{P}_1(\cdot)}\sum_{i\in\calX}\sum_{j\in\calX}\ell_j\Xi_{ij}\mu_i .
\end{equation}
In addition,
\begin{equation}\label{eq.posnegvar}
\sum_{i\in\calX}\sum_{j\in\calX}\ell_j\Xi_{ij}\mu_i=\sum_{i\in\calX}\sum_{j\in\calX}\ell_j\Xi_{ij}^+\mu_i-\sum_{i\in\calX}\sum_{j\in\calX}\ell_j\Xi_{ij}^-\mu_i.
\end{equation}

The solution of Problem \ref{problem3} is obtained by identifying the partition of $\calX$ into disjoint sets $\{\calX^0,\calX_0,\calX_1,\dots,\calX_k\}$ and the transitions on this partition. The main idea is to express $\Xi_{i\bullet}$ as the difference of its positive and negative part and then find upper and lower bounds on the transition probabilities of $\calX^0$ and $\calX\setminus\calX^0$ which are achievable. Closed form expressions of the transition probability measures, on these sets, which achieve the bounds are derived. 

Note that, if we replace the maximization in \eqref{MCproblem1} with minimization, then the solution of the new problem is obtained precisely as that of Problem \ref{problem3}, but with a reverse computation of the partition of the space $\calX$ and the mass of the transition probability on the partition moving in the opposite direction. 

The following Theorem characterizes the solution of Problem \ref{problem3}.

\begin{theorem}\label{MCthmocc}
The solution of Problem \ref{problem3} is given by
\begin{equation}
\sum_{i\in \calX}\sum_{j\in \calX}\ell_j \Phi^{\dagger}_{ij}\mu_i=\ell_{\max}\sum_{i\in \calX^0}\sum_{j\in \calX}\mu_j\Phi^{\dagger}_{ji}+\ell_{\min}\sum_{i\in \calX_0}\sum_{j\in \calX}\mu_j\Phi^{\dagger}_{ji}
+\sum_{k=1}^r\ell(\calX_k)\sum_{i\in \calX_k}\sum_{j\in\calX}\mu_j\Phi^{\dagger}_{ji} ,
\end{equation}
where for any $i\in\calX$,
\begin{subequations}\label{eq.55}
\begin{align}
\Phi^{\dagger}_{ij}&=P_{ij}+\frac{\alpha_i}{2|\calX^0|}, \quad \forall j\in\calX^0,\\
\Phi^{\dagger}_{ij}&=\Big(P_{ij}-\frac{\alpha_i}{2|\calX_0|}\Big)^+,\quad \forall j\in\calX_0,\\
\Phi^{\dagger}_{ij}&=\Big(P_{ij}-\Big(\frac{\alpha_i}{2|\calX_k|}-\sum_{j=1}^k\sum_{z\in\calX_{j-1}}P_{iz}\Big)^+\Big)^+, \quad \forall j\in\calX_k \\
 \alpha_i&=\min(R,R_{\max,i}),\quad R_{\max,i}=2(1-\sum_{j\in\calX^0}P_{ij}) \label{eq.55d},
\end{align}
\end{subequations}
$k=1,2,\dots,r$ and $r$ is the number of $\calX_k$ sets which is at most $|\calX\setminus \calX^0\cup\calX_0|$.
Once the $\Phi^{\dagger}$ matrix is constructed as a function of TV parameter $R$, then the transition matrix $\Phi$ which solves \eqref{MCproblem1} is given by removing all zero columns and the respective rows of $\Phi^{\dagger}$ matrix.
\end{theorem}

\begin{IEEEproof}
See the Appendix. 
\end{IEEEproof}
 
Clearly, the optimal transition matrix $\Phi$ is obtained via a \textit{water-filling} solution. 

%
%
%
%
\section{{Method 2}: Solution of approximation problems}\label{sec.sol.opt.problem.}
In this section, we recall some results from \cite{ctlthem2013}, which are vital in providing the solution of Problem \eqref{max.aver.equation}, and consequently the solution of approximation Problems \eqref{problem.occ} and \eqref{problem.entr}. 

First recall, from Section \ref{sec.sol.opt.problem.2}, the definitions of the support sets $\calX^0$, $\calX_0$, $\calX_k$ and the definitions of the corresponding values of the sequence on these sets given by $\ell_{\max}$, $\ell_{\min}$ and $\ell(\calX_k)$.

Given $\ell\in \mathbb{R}_+^{|\calX|}$, $\mu\in {\mathbb P}({\calX})$, it is shown in \cite{ctlthem2013}, that the solution of optimization \eqref{max.aver.equation} is given by
\begin{equation} 
\mathbb{L}(\nu^*)=\ell_{\max}\nu^*(\calX^0)+\ell_{\min}\nu^*(\calX_0)+\sum_{k=1}^r\ell(\calX_k)\nu^*(\calX_k),\label{mp}
\end{equation}
\noi and the optimal probabilities are obtained via \textit{water-filling}, as follows 
\begin{subequations}\label{all3}
\begin{align}{l}
 \nu^*(\calX^0)&\triangleq \sum_{i\in\calX^0}\nu_i^*=\sum_{i\in\calX^0}\mu_i+\frac{\alpha}{2},\label{all3a}\\
 \nu^*(\calX_0)&\triangleq \sum_{i\in\calX_0}\nu_i^*=\Big(\sum_{i\in\calX_0}\mu_i-\frac{\alpha}{2}\Big)^+,\label{all3b}\\
 \nu^*(\calX_k)&\triangleq \sum_{i\in\calX_k} \nu_i^*=\Big(\sum_{i\in\calX_k} \mu_i-\Big(\frac{\alpha}{2}-\sum_{j=1}^k\sum_{i\in\calX_{j-1}}\mu_i\Big)^{+}\Big)^+, \label{all3c}\\
\alpha&=\min\left(R,R_{\max}\right),\qquad R_{\max}\tri 2(1-\sum_{i\in \calX^0}\mu_i),\label{all3d}
\end{align}
\end{subequations}
where, $k=1,2,\hdots,r$ and $r$ is the number of $\calX_k$ sets which is at most $|\calX\setminus\calX^0\cup\calX_0|$. The optimal probabilities given by \eqref{all3a}-\eqref{all3c}, can be expressed in matrix form as follows
\bea \label{MC2}
\nu^*=\mu Q^{\dagger}=\mu PQ^{\dagger}.
\eea
 In Sections \ref{subsec.occ.problem.} and \ref{subsec.entr.problem.}, we provide algorithms for constructing the desired $Q^{\dagger}$ matrix for the optimizations \eqref{max.occupancy} and \eqref{max.entropy}, respectively.
 
 \begin{remark}\label{remar11}
 The identification of the support sets $\calX^0$, $\calX_0$ and $\calX_k$, $k=1,2,\dots,r$, is based on the values of $\ell_i$'s, $\forall i\in\calX$. If the cardinality of any of the support sets is greater than one, i.e., $|\calX^0|>1$, and $\ell_i=\ell_{i+1}=\dots$, $\forall i,i+1,\dots\in\calX^0$ then by \eqref{all3a} 
\begin{subequations}\label{eq:all16}
  \begin{align}
  \nu_i^*&=\frac{\nu^*(\calX^0)}{|\calX^0|}, \quad\forall i\in\calX^0,\\
 \intertext{and similarly for the rest, that is, if $|\calX_0|>1$ then}
  \nu_i^*&=\frac{\nu^*(\calX_0)}{|\calX_0|}, \quad\forall i\in\calX_0,\\
   \intertext{and if $|\calX_k|>1$, for $k=1,\dots,r$, then}
   \nu_i^*&=\frac{\nu^*(\calX_k)}{|\calX_k|}, \quad\forall i\in\calX_k.
  \end{align}
\end{subequations}   
The resulting optimal probability $\nu^*$ is a $(2+r)$ row vector and hence, by \eqref{MC2} $Q^\dagger$ is an $|\calX|\times(2+r)$ matrix. Then by employing \eqref{eq:all16} we extract the optimal probabilities $\nu_i^*$ for all $i\in\calX$, which are then used in definition of the optimal partition functions (see Definition \ref{part.function.occ} and \ref{part.function.entr}).

For the approximation based on occupancy distribution, we let the matrix $Q^\dagger$ to be an $|\calX|\times|\calX|$ matrix, instead of an $|\calX|\times (2+r)$ matrix. The reason for doing so, is that we want to take into account the cases for which $\ell_i$'s, $\forall i\in\calX$, might be defined to represent a cost or profit etc., whenever a node is visited. In such cases, \eqref{eq:all16} is not valid anymore, since $\ell_i=\ell_j$ does not necessarily imply $\mu_i=\mu_j$, $\forall i,j\in\calX$. As we will show in Section \ref{subsec.occ.problem.}, Algorithm \ref{algorithm1} constructs a $Q^\dagger$ matrix which in addition to occupancy distribution, considers those alternative cases as well.
 \end{remark}

 By Definition \ref{lower.dim.process1} and \ref{lower.dim.process2}, the approximated probability vector $\bar{\nu}\in \mathbb{P}(\calY)$ is readily available and satisfies \bea 
\bar{\nu}=\mu Q=\mu PQ , 
\eea 
where $Q$ matrix is modified accordingly.

\begin{figure}[!h]
\centering
\includegraphics[width=.8\linewidth]{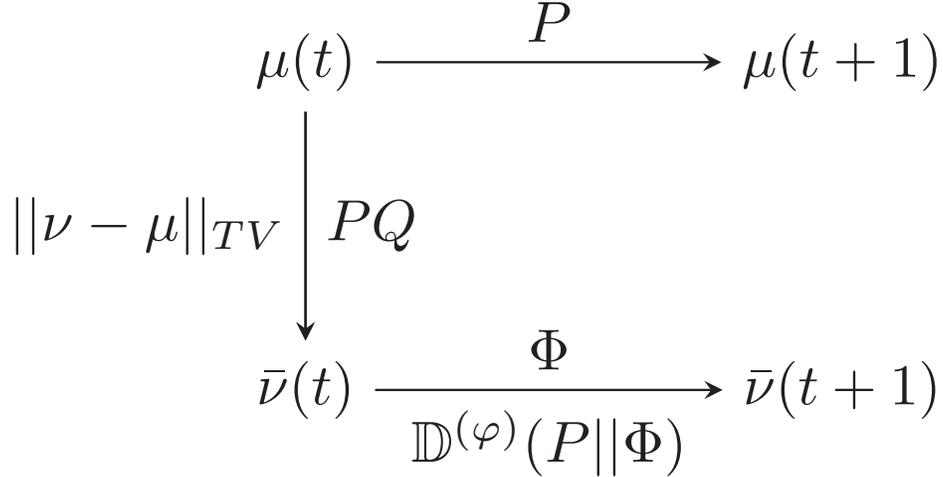}
\caption[]{\emph{Method 2}. Diagram that shows the relationship of the initial and the lower probability distributions.}\label{diag1}
\end{figure}

Once the reduced state process is obtained, we utilize its solution to solve the optimizations \eqref{KLoccupancy} and \eqref{KLentropy}. The relation between $\mu(t),\mu(t+1)\in\mathbb{P}(\calX)$ and $\bar{\nu}(t), \bar{\nu}(t+1)\in \mathbb{P}(\calY)$ is shown in Fig.\ref{diag1}.

\subsection{Solution of approximation problem based on occupancy distribution}\label{subsec.occ.problem.}
In this section, we first give an algorithm to construct the $Q^{\dagger}$ matrix which solves \eqref{max.occupancy}. Then, under an additional assumption that the reduced process is also Markov, we give the solution of \eqref{KLoccupancy}.
 
Let $k=0,1,\hdots,r-1$, where $r$ denotes the number of $\calX_k$ sets, that is, $1\leq r \leq |\calX\setminus \calX^0|$ (note that, $\calX_0$ set is included). For all $j=1,2,\dots,|\calX_k|$, $\calX_{k,j}\triangleq \{jth \ \mbox{element of}\ \calX_k \ \mbox{set}\}$, (note that, if $|\calX_k|=1$ then $\calX_{k,j}=\calX_k$). Similarly, $\calX^{0,j}\triangleq \{jth \ \mbox{element of}\ \calX^0\ \mbox{set}\}$, (note that, if $|\calX^0|=1$ then $\calX^{0,j}=\calX^0$).

\begin{algorithm}\label{algorithm1}
\noi
 \ben
 \item Initialization step:
\ben
\item Arrange $\ell_i$, $i\in\calX$, in a descending order.
\item Identify the support sets $\calX^0$, $\calX_0$ and $\calX_k$ for all $k\in\{1,2,\dots,|\calX\setminus\calX^0\cup\calX_0|\}$.
\item Calculate the value of $r$.
\een
\een

For any $R\in[0,2]$:
\ben 
\item[2)] Step.1 (Indicator functions): 
\ben
\item Let 
\een
\een
\begin{equation*} 
\mu^R(\calX^0)\triangleq \sum\limits_{i\in\calX^0}\mu_i+\frac{R}{2}.
\end{equation*}
\ben
\item[] \ben \item[] Define \een
\een
\begin{equation}\label{alg1ind1}
  I^{\calX^0} \triangleq  \left\{ \,
\begin{IEEEeqnarraybox}[][c]{l?s}
\IEEEstrut
1, & if $\mu^R\left(\calX^0\right)\geq 1$, \\
0 , & otherwise.
\IEEEstrut
\end{IEEEeqnarraybox}
\right.
\end{equation}
\ben
\item[]
\ben
\item[b)] For $k=0,1,\dots,r-1$ let
\een
\een
\begin{equation*}
\mu^{R}(\calX_k)\triangleq\sum\limits_{ \mathclap{j=0}}^k\sum\limits_{i\in\calX_j}{\mu_i}-\frac{R}{2}.
\end{equation*}
\ben
\item[] \ben 
\item[] Define 
\een
\een
\begin{equation*}
I^{\calX_k} \triangleq \left\{ \,
\begin{IEEEeqnarraybox}[][c]{l?s}
\IEEEstrut
1, & if $\mu^R\left(\calX_k\right) \geq 0$,  \\
0, & otherwise.
\IEEEstrut
\end{IEEEeqnarraybox}
\right.
\end{equation*}
\begin{equation*}
I^{\calX_{[0,k-1]}} \triangleq \left\{ \,
\begin{IEEEeqnarraybox}[][c]{l?s}
\IEEEstrut
1, & if $\mu^R\left(\calX_i\right) <0, \ \forall i=0,1,\dots,k-1$  \\
0, & otherwise,
\IEEEstrut
\end{IEEEeqnarraybox}
\right.
\end{equation*}
\ben
\item[] \ben 
\item[] and
\een
\een
\begin{equation}\label{alg1ind2}
I^{\calX_k,\calX_{[0,k-1]}}=I^{\calX_k}I^{\calX_{[0,k-1]}}.
\end{equation}
\ben
\item[]
\ben
\item[c)] For $k=0,1,\dots,r-1$, if $|\calX_k|>1$, then for all $j{=}1\dots,|\calX_k|$, let 
\een
\een
\begin{IEEEeqnarray*}{l}
\mu^R(\calX_{k,j})\triangleq \mu_{\calX_{k,j}}-\frac{(R/2-\sum_{i\in\cup_{j=0}^{k-1}\calX_j}\mu_i)}{|\calX_k|}.
\end{IEEEeqnarray*}
\ben
\item[] \ben 
\item[] Define 
\een
\een
\begin{equation}\label{alg1ind3}
I^{\calX_{k,j}} \triangleq \left\{ \,
\begin{IEEEeqnarraybox}[][c]{l?s}
\IEEEstrut
1, & if $\mu^R\left(\calX_{k,j}\right)\geq 0$, \\
0 , & otherwise,
\IEEEstrut
\end{IEEEeqnarraybox}
\right.
\end{equation}
\ben 
\item[3)] Step.2 (The $Q^\dagger$ matrix): 

Let $Q^\dagger$ be an $|\calX|\times|\calX|$ matrix and $i=1,2,\hdots,|\calX|$ to denote the $i$th column of $Q^\dagger$ matrix.
\ben
\item[a)]  For all $i\in \calX^0$, the elements of the $ith$ column are given as follows.
\een
\een
\ben
\item[] \ben
\item[]
\ben
\item[i)] Let the $(Q^\dagger)_{i,i}$ element be equal to
\een
\een
\een
\begin{equation} \sum_{k=0}^{r-1}I^{\calX_k,\calX_{[0,k-1]}}\Big(1+\frac{R/2}{|\calX^0|}\Big)+
I^{\calX^0}\frac{(\mu_{\calX^{0,i}}+\frac{\sum_{j\in\calX\setminus\calX^0}\mu_j}{|\calX^0|})}{\mu_{\calX^{0,i}}}.
\end{equation}
\ben
\item[] \ben
\item[]
\ben
\item[ii)] Let all the remaining elements of the $i$th column be equal to
\een
\een
\een
 \begin{equation}
\sum_{k=0}^{r-1}I^{\calX_k,\calX_{[0,k-1]}}\frac{R/2}{|\calX^0|}.
\end{equation}
\ben \item[] \ben
\item[b)] For all $i\in\calX_k$, $k=0,1,\dots,r-1$, and $j\in\Big\{\psi\in\{1,2,\dots,|\calX_k|\}: i\in\calX_k \ \mbox{is in the}\\ \ \psi th \ \mbox{position on}\  \calX_k \ \mbox{set}\Big\}$, the elements of the $ith$ column are as follows.
\een \een

\ben
\item[] \ben
\item[]
\ben
\item[i)] Let the $(Q^\dagger)_{i,i}$ element be equal to
\een
\een
\een
\begin{equation}
\sum_{j=0}^{k-1}I^{\calX_j,\calX_{[0,j-1]}}+I^{\calX_{k,j}}I^{\calX_k,\calX_{[0,k-1]}}\Big(1-\frac{R/2}{\sum_{j=1}^{|\calX_k|}I^{\calX_{k,j}}}\Big).
\end{equation}
\ben
\item[] \ben
\item[]
\ben
\item[ii)] If $|\calX_k|>1$, then for all $z\in\calX_k\setminus\calX_{k,j}$, let the $(Q^\dagger)_{z,i}$ element be equal to
\een
\een
\een
\begin{equation}
I^{\calX_{k,j}}I^{\calX_k,\calX_{[0,k-1]}}\Big\{
\prod_{j=1}^{|\calX_k|}I^{\calX_{k,j}} \Big(\frac{-R/2}{\sum_{j=1}^{|\calX_k|}I^{\calX_{k,j}}}\Big)
+\Big(1-\frac{R/2}{\sum_{j=1}^{|\calX_k|}I^{\calX_{k,j}}}\Big)\Big(1-\prod_{j=1}^{|\calX_k|}I^{\calX_{k,j}}\Big)\Big\} .
\end{equation}
\ben
\item[] \ben
\item[]
\ben
\item[iii)] For all $z\in\calX\setminus\calX^0\cup\calX_k$ and only if $z>i$ let the $(Q^\dagger)_{z,i}$ element be equal to
\een
\een
\een
\begin{equation}
I^{\calX_{k,j}}I^{\calX_k,\calX_{[0,k-1]}}\Big\{
\prod_{j=1}^{|\calX_k|}I^{\calX_{k,j}}\Big(\frac{1}{|\calX_k|}-\frac{R/2}{\sum_{j=1}^{|\calX_k|}I^{\calX_{k,j}}}\Big)\\
{+}\Big(1{-}\frac{R/2}{\sum_{j=1}^{|\calX_k|}I^{\calX_{k,j}}}\Big)\Big(1{-}\prod_{j=1}^{|\calX_k|}I^{\calX_{k,j}}\Big)\Big\}.
\end{equation}
\ben
\item[] \ben
\item[]
\ben
\item[iv)] Let all the remaining elements of the $i$th column be equal to
\een
\een
\een
\begin{equation}
I^{\calX_{k,j}}I^{\calX_k,\calX_{[0,k-1]}}\Big(\frac{-R/2}{\sum_{j=1}^{|\calX_k|}I^{\calX_{k,j}}}\Big) .
\end{equation}
\end{algorithm}

Once the $Q^{\dagger}$ matrix is constructed, as a function of the TV parameter $R$, then by \eqref{MC2} the resulting optimal probability, $\nu^*$, is an $1\times|\calX|$ row vector. However, recall from Remark \ref{remar11} that by definition $\nu^*$ is just an $1\times(2+r)$ row vector. By using all the information that the support sets provide to us we can easily transform the $1\times|\calX|$ row vector to an $1\times(2+r)$ row vector, by simply adding together the optimal probabilities, $\nu^*_i$, $\forall i\in\calX$, which belong to the same support sets. Given the optimal solution of optimization \eqref{max.occupancy}, then by Definition \ref{lower.dim.process1} the lower dimensional process $\{Y_t:t=0,1\dots\}$ with invariant distribution $\bar{\nu}$ is obtained, either by removing all zero elements of $\nu^*\in \mathbb{P}(\calX)$, or by defining a $Q$ matrix to be equal to $Q^{\dagger}$ after the deletion of all zero columns, and hence \begin{equation}
\bar{\nu}=\mu Q=\mu PQ ,
\end{equation}
where the dimensions of $Q$ matrix are based on the value of the TV parameter $R\in[0,2]$.

Before we proceed with the solution of \eqref{KLoccupancy}, we provide a simple, yet useful example in order to explain each step of Algorithm \ref{algorithm1}.

\begin{example}\label{ex.alg.occupancy}
 Let $\ell=[\ell_1 \ \ell_2\  \ell_3\ \ell_4]$, where $\ell_1>\ell_2>\ell_3>\ell_4$, and $|\calX|=4$. For simplicity it is assumed that the optimum probabilities $\nu_i^*$, $i\in \calX$, as a function of $R$ are known, as presented in Fig.\ref{fig1000}. 
 
Initialization step. The support sets are equal to $\calX^0=\{1\}$, $\calX_0=\{4\}$, $\calX_1=\{3\}$ and $\calX_2=\{2\}$. The number of $\calX_k$ sets is equal to $r=3$. 

 Step.1 From \eqref{alg1ind1}, the indicator function $I^{\calX^0}$ is given by
\begin{equation*}
I^{\calX^0}\triangleq \left\{ \,
\begin{IEEEeqnarraybox}[][c]{l?s}
\IEEEstrut
1, & if $\mu_1+\frac{R}{2}\geq 1$,\\
0, & otherwise.
\IEEEstrut
\end{IEEEeqnarraybox}
\right.
\end{equation*}
From \eqref{alg1ind2}, the indicator functions $I^{\calX_0}$, $I^{\calX_1,\calX_{[0,0]}}$ and $I^{\calX_2,\calX_{[0,1]}}$ are given by
\begin{equation*}
I^{\calX_0} \triangleq \left\{ \,
\begin{IEEEeqnarraybox}[][c]{l?s}
\IEEEstrut
1, & if $\mu_4-\frac{R}{2}\geq 0$,\\
0, & otherwise.
\IEEEstrut
\end{IEEEeqnarraybox}
\right.
 \end{equation*}
\begin{equation*}
I^{\calX_1,\calX_{[0,0]}}\triangleq \left\{ \,
\begin{IEEEeqnarraybox}[][c]{l?s}
\IEEEstrut
1, & if $\mu_3+\mu_4-\frac{R}{2}\geq 0$ and $\mu_4-\frac{R}{2}\leq 0$,\\
0, & otherwise.
\IEEEstrut
\end{IEEEeqnarraybox}
\right.
\end{equation*}
\begin{equation*}
I^{\calX_2,\calX_{[0,1]}}\triangleq \left\{ \,
\begin{IEEEeqnarraybox}[][c]{l?s}
\IEEEstrut
1, & if $\mu_2+\mu_3+\mu_4-\frac{R}{2}\geq 0$ \\
& and $\mu_3{+}\mu_4{-}\frac{R}{2}\leq 0$ and $\mu_4{-}\frac{R}{2}\leq 0$,\\
0, & otherwise.
\IEEEstrut
\end{IEEEeqnarraybox}
\right.
\end{equation*}
The values of the indicator functions for $R\in [0,2]$ are given below.\\
\begin{minipage}[t]{0.24\linewidth}%
\centering
\underline{$0\leq R< R_1$}
\begin{IEEEeqnarray*}{rCl}
I^{\mathcal{X}^0}&{=}&0\\
I^{\mathcal{X}_0}&{=}&1\\
I^{\mathcal{X}_1,\mathcal{X}_{[0,0]}}&{=}&0\\
I^{\mathcal{X}_2,\mathcal{X}_{[0,1]}}&{=}&0
\end{IEEEeqnarray*}
\end{minipage}
\begin{minipage}[t]{0.24\linewidth}%
\centering
\underline{$R_1\leq R< R_2$}
\begin{IEEEeqnarray*}{rCl}
I^{\mathcal{X}^0}&{=}&0\\
I^{\mathcal{X}_0}&{=}&0\\
I^{\mathcal{X}_1,\mathcal{X}_{[0,0]}}&{=}&1\\
I^{\mathcal{X}_2,\mathcal{X}_{[0,1]}}&{=}&0
\end{IEEEeqnarray*}
\end{minipage}
\begin{minipage}[t]{0.24\linewidth}%
\centering
\underline{$R_2\leq R<R_3$}
\begin{IEEEeqnarray*}{rCl}
I^{\mathcal{X}^0}&{=}&0\\
I^{\mathcal{X}_0}&{=}&0\\
I^{\mathcal{X}_1,\mathcal{X}_{[0,0]}}&{=}&0\\
I^{\mathcal{X}_2,\mathcal{X}_{[0,1]}}&{=}&1
\end{IEEEeqnarray*}
\end{minipage}
\begin{minipage}[t]{0.24\linewidth}%
\centering
\underline{$ R_3\leq R\leq 2$}
\begin{IEEEeqnarray*}{rCl}
I^{\mathcal{X}^0}&{=}&1\\
I^{\mathcal{X}_0}&{=}&0\\
I^{\mathcal{X}_1,\mathcal{X}_{[0,0]}}&{=}&0\\
I^{\mathcal{X}_2,\mathcal{X}_{[0,1]}}&{=}&0\\
\end{IEEEeqnarray*}
\end{minipage}

\noi For $0 \leq R< R_1$, all indicator functions are equal to one, except the one which corresponds to $\calX_0$ set, that is, $I^{\calX_0}=1$. As soon as $\mu^R(\calX_0)=0$, then $I^{\calX_0}$ becomes equal to zero and $I^{\calX_1,\calX_{[0,0]}}$ equal to one. This procedure is repeated until the value of $R=R_{\max}=R_3$, see Fig.\ref{fig1000}, in which $I^{\calX^0}$ becomes equal to one, and all other indicator functions equal to zero, and $I^{\calX^0}$ remains active for all $R\geq R_{\max}=R_3$.

Step.2 Let $Q^{\dagger}$ be an $4\times 4$ matrix. For $0 \leq R< R_1$,
\bes Q^{\dagger}=\left(
         \begin{array}{cccc}
           1+R/2 & 0 & 0 & -R/2 \\
           R/2 & 1 & 0 & -R/2 \\
          R/2 & 0 & 1 &-R/2 \\
           R/2 & 0 & 0 &1-R/2 \\
         \end{array}
       \right),\ees
and since no zero column exist then $Q^{\dagger}=Q$. For $ R_1\leq R< R_2$,
\bes Q^{\dagger}=\left(
         \begin{array}{cccc}
           1+R/2 & 0  & -R/2 & 0\\
          R/2 & 1 & -R/2 & 0 \\
          R/2 & 0  &1-R/2 & 0\\
          R/2 & 0  &1-R/2 & 0\\
         \end{array}
       \right)\Longrightarrow Q=\left(
         \begin{array}{cccc}
           1+R/2 &0 & -R/2 \\
           R/2 &1 & -R/2 \\
          R/2 &0 &1-R/2 \\
           R/2 &0 &1-R/2 \\
         \end{array}
       \right).\ees   
For $ R_2\leq R< R_3$,
\bes Q^{\dagger}=\left(
         \begin{array}{cccc}
           1{+}R/2   & {-}R/2& 0 & 0\\
          R/2  & 1{-}R/2& 0 & 0 \\
          R/2   &1{-}R/2& 0 & 0\\
           R/2   &1{-}R/2& 0 & 0\\
         \end{array}
       \right)\Longrightarrow Q=\left(
         \begin{array}{cccc}
           1{+}R/2   & {-}R/2\\
          R/2  & 1{-}R/2 \\
          R/2   &1{-}R/2\\
           R/2   &1{-}R/2\\
         \end{array}
       \right).\ees
For $ R\geq R_3$,
\bes Q^{\dagger}=\left(
         \begin{array}{cccc}
           \frac{1}{\mu_1}   & 0& 0 & 0\\
          0  & 0& 0 & 0 \\
          0   &0& 0 & 0\\
           0  &0& 0 & 0\\
         \end{array}
       \right)\Longrightarrow Q=\left(
         \begin{array}{cccc}
          \frac{1}{\mu_1}  \\
           0   \\
          0  \\
           0  \\
         \end{array}
       \right).\ees      
\begin{figure}[h!]
\centering
\includegraphics{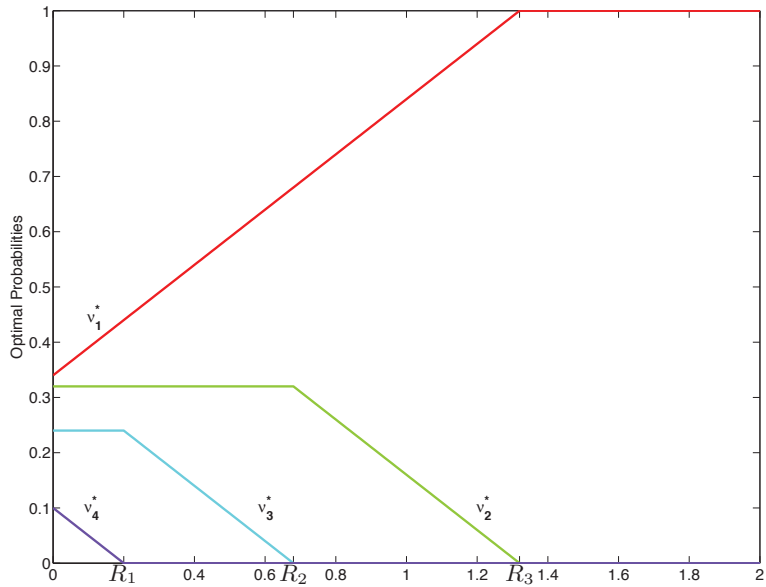}
\caption[]{Optimal probabilities as a function of $R$.}
\label{fig1000}
\end{figure}
Note that, the number of columns of $Q$ matrix is based on the value of total variation parameter $R$. For $0\leq R <R_1$, its dimension is equal to $(|\calX|)\times (1+r)$. Whenever an indicator function becomes equal to zero, all elements of the respective column become equal to zero, and hence the column is deleted, until $R\geq R_3$, where the $Q$ matrix will be transformed into a column vector of dimension $(|\calX|)\times (1)$. 
\end{example}

 Next, we proceed with the solution of \eqref{KLoccupancy}, by letting $\bar{\nu}\in\mathbb{P}(\calY)$ to denote the invariant distribution of a lower dimensional Markov process $(\bar{\nu},\Phi)$. As mentioned in \cite{Deng:11}, the main difficulty in solving \eqref{KLoccupancy} is in finding an optimal partition function $\varphi$. However, once an optimal partition is given then the solution of $\Phi$ can be easily obtained. Toward this end, next we define an optimal partition function  for the approximation problem based on occupancy distribution at values of TV parameter $R$ for which a reduction of the states occurs (i.e., see Example \ref{ex.alg.occupancy}, Fig.\ref{fig1000}, for values of $R=R_1, R_2$ and $R_3$).
\begin{definition}\label{part.function.occ}(Partition function) Let $\calX$ and $\calY$ be two finite dimensional state-spaces with $|\calY|<|\calX|$. Define a surjective (partition) function $\varphi:\calX\longmapsto \calY$ as follows. \begin{align*}
\forall i\in \calX^0,\quad \varphi(i)&=1\in \calY ,\\ 
\forall i\in \calX\setminus\calX^0,\quad  \varphi(i) &= \left\{
  \begin{array}{l l}
    1, & \ \text{if $\nu^*_i=0$,}\\
    k\in \calY , & \ \text{if $\nu^*_i>0$.}
  \end{array} \right.
\end{align*}
\end{definition}

Note that, once the optimal probabilities $\nu_i^*$, $\forall i\in\calX$ are obtained, we can easily identify the values of $R$ for which a reduction of the states occurs. In addition, since the solution behavior of \eqref{max.occupancy} is to remove probability mass from states with the smallest invariant probability and strengthening the states with the highest invariant probability, this property of the partition function $\varphi$ is intuitive and expected.

Next, we reproduce the main theorem of \cite{Meyn09}, which gives the solution of $\Phi$ that solves \eqref{KLoccupancy}.
\begin{theorem}
Let $(\mu,P,\calX)$ be a given FSM process and $\varphi$ be the partition function of Definition \ref{part.function.occ}. For optimization \eqref{KLoccupancy}, the solution of $\Phi$ is given by
\begin{equation}
\Phi_{kl}=\frac{u^{(k)}\Pi Pu^{(\ell)'}}{\bar{\nu}_k},\quad k,\ell\in \calY ,
\end{equation}
where $\Pi=diag(\mu)$, $u^{(k)'}$ is the transpose of $u^{(k)}$, and $u^{(k)}$ is a $1\times|\calX|$ row vector defined by
\begin{equation} u_i^{(k)} = \left\{
  \begin{array}{l l}
    1, & \quad \text{if $\varphi(i)=k$,}\\
   0, & \quad \text{otherwise.}
  \end{array} \right.\end{equation}
\end{theorem}

\begin{IEEEproof}
See \cite{Deng:11}.
\end{IEEEproof}

\subsection{Solution of Approximation problem based on maximum entropy principle}\label{subsec.entr.problem.}
In this subsection, we first give an algorithm to construct the $Q$ matrix which solves  \eqref{max.entropy}. Then, under the assumption that the reduced process is also Markov, we give the solution of \eqref{KLentropy}. Before giving the algorithm, we introduce some notation.

Let $r$ denote the number of $\calX_k$ sets, that is, $1\leq r\leq |\calX\setminus \calX^0\cup \calX_0|$ (note that, $\calX_0$ set is excluded, in contrast with the definition of $r$ in Section \ref{subsec.occ.problem.}). Furthermore, let $r^+$ and $r^-$ denote the number of $\mu_i$, $i\in\calX$, such that $\mu_i\geq \frac{1}{|\calX|}$ and $\mu_i< \frac{1}{|\calX|}$, respectively. In addition, $\mu_i\neq \mu_j$ should also be satisfied for all $i\neq j$, $i,j\in \calX$.

\begin{remark}\label{remalg2}
The initialization step of the following algorithm is performed by letting $R=0$. In this case, $\nu_i=\mu_i$, $\forall i\in \calX$, and hence, $\ell_i\triangleq -\log \nu_i=-\log\mu_i$.
\end{remark}
\begin{algorithm}\label{algorithm2}
\noi
 \ben
 \item Initialization step:
\ben
\item  Arrange $\mu_i$, $i\in\calX$, in a descending order and let $R=0$.
\item Identify the support sets $\calX^0$, $\calX_0$ and $\calX_k$ for all $k\in\{1,2,\dots,|\calX\setminus\calX^0\cup\calX_0|\}$.
\item Calculate the value of $r$, $r^-$ and $r^+$.
\een
\een

For any $R\in[0,2]$:
\ben 
\item[2)] Step.1 (Indicator functions): 
\ben
\item For $k=1,2\hdots,r^-{-}1$ let
\een
\een
\begin{equation*}
\mu^R_-(\calX_k)\triangleq \frac{\sum_{i\in\cup_{j=0}^{k-1}\calX_j}\mu_i-R/2}{\sum_{j=0}^{k-1}|\calX_j|}.
\end{equation*}
\ben
\item[] \ben 
\item[] Define
\een
\een
\begin{equation}\label{alg2ind1}
I^{\calX_k}_{-}\triangleq \left\{ \,
\begin{IEEEeqnarraybox}[][c]{l?s}
\IEEEstrut
1, & if $\mu^R_-(\calX_k)\leq \frac{\sum_{i\in\calX_{k}}\mu_i}{|\calX_k|}$,\\
0, & otherwise.
\IEEEstrut
\end{IEEEeqnarraybox}
\right.
\end{equation}
\ben \item[] \ben \item[]
For $k=r^-$ let \een \een
\begin{equation*}
\mu^R_-(\calX_{r^-})\triangleq \frac{\sum_{i\in\cup_{j=0}^{r^-{-}1}\calX_j}\mu_i-R/2}{\sum_{j=0}^{r^-{-}1}|\calX_j|}.
\end{equation*}
\ben
\item[] \ben 
\item[] Define
\een
\een
\begin{equation}\label{alg2ind2}
I^{\calX_{r^-}}_{-}\triangleq \left\{ \,
\begin{IEEEeqnarraybox}[][c]{l?s}
\IEEEstrut
1, & if $\mu^R_-(\calX_{r^-})\leq\frac{1}{|\calX|}$,\\
0, & otherwise.
\IEEEstrut
\end{IEEEeqnarraybox}
\right.
\end{equation}
\ben \item[] \ben \item[b)] For $k=1,2\hdots,r^+{-}1$ let\een \een
\begin{equation*}
\mu^R_+(\calX_k)\triangleq \frac{\sum_{i\in\calX\setminus\cup_{j=r}^{k-1}\calX_{r-j}}\mu_i+R/2}{|\calX\setminus\cup_{j=r}^{k-1}\calX_{r-j}|}.
\end{equation*}
\ben
\item[] \ben 
\item[] Define
\een
\een
\begin{equation}\label{alg2ind3}
 I^{\calX_k}_{+}\triangleq  \left\{ \,
\begin{IEEEeqnarraybox}[][c]{l?s}
\IEEEstrut
1, & if $\mu_+^R(\calX_k)\geq \frac{\sum_{i\in\calX_{r-k+1}}\mu_i}{|\calX_{r-k+1}|}$, \\
0, & otherwise.
\IEEEstrut
\end{IEEEeqnarraybox}
\right.
\end{equation}
\ben \item[] \ben \item[] For $k=r^+$ let \een \een
\begin{equation*}
\mu^R_+(\calX_{r^+})\triangleq \displaystyle\frac{\sum_{i\in\calX\setminus\cup_{j=r}^{{r^+}{-}1}\calX_{r-j}}\mu_i+\frac{R}{2}}{|\calX\setminus\cup_{j=r}^{{r^+}{-}1}\calX_{r-j}|}.
\end{equation*}
\ben
\item[] \ben 
\item[] Define
\een
\een
\begin{equation}\label{alg2ind4}
 I^{\calX_{r^+}}_{+}\triangleq  \left\{ \,
\begin{IEEEeqnarraybox}[][c]{l?s}
\IEEEstrut
1, & if $\mu_+^R(\calX_{r^+})\geq \frac{1}{|\calX|}$, \\
0, & otherwise.
\IEEEstrut
\end{IEEEeqnarraybox}
\right.
\end{equation}

\ben 
\item[3)] Step.2 (The $Q^\dagger$ matrix): 

Let $Q^{\dagger}$ be an $(|\calX|)\times (2+r)$ matrix.

\ben \item The elements of the first column are given as follows.

\ben
\item[i)] For all $i\in\calX_0$, let the $(Q^{\dagger})_{i,1}$ be equal to 
\een
\een
\een
\begin{equation}
\frac{1-R/2}{|\calX_0|+\sum_{j=1}^{r^{\downarrow}-1}I^{\calX_j}_{-}|\calX_j|}\Big(I^{\calX_{r^-}}_{-}\Big)^\mathsf{c}+\frac{I^{\calX_{r^-}}_{-}}{|\calX|}.
\end{equation}

\ben \item[] \ben \item[] \ben \item[ii)]For all $i\in\calX_k$, $k=1,2,\hdots,r^-{-}1$, let the $(Q^{\dagger})_{i,1}$ be equal to \een \een \een
\begin{equation}
\frac{I^{\calX_k}_{-}-R/2}{|\calX_0|+\sum_{j=1}^{r^{\downarrow}-1}I^{\calX_j}_{-}|\calX_j|}\Big(I^{\calX_{r^-}}_{-}\Big)^\mathsf{c}+\frac{I^{\calX_{r^-}}_{-}}{|\calX|}.
\end{equation}
\ben \item[] \ben \item[] \ben \item[iii)] Let all the remaining elements be equal to \een \een \een
\begin{equation} \frac{-R/2}{|\calX_0|+\sum_{j=1}^{r^{\downarrow}-1}I^{\calX_j}_{-}|\calX_j|}\Big(I^{\calX_{r^-}}_{-}\Big)^\mathsf{c}+\frac{I^{\calX_{r^-}}_{-}}{|\calX|}.\end{equation}

\ben \item[] \ben \item[b)] The elements of the last column are given by 
\ben \item[i)] For all $i\in\calX^0$, let the $(Q^{\dagger})_{i,r+2}$ be equal to \een
\een \een
\begin{equation}
\frac{1+R/2}{|\calX^0|+\sum_{j=1}^{r^{\uparrow}-1} I^{\calX_j}_{+}|\calX_{r-j+1}|}(I^{\calX_{r^+}}_{+})^\mathsf{c}.
\end{equation}
\ben \item[]\ben  \item[]\ben \item[ii)] For all $i\in\calX_{r-k+1}$, $k=1,2,\hdots,r^{\uparrow}-1$ let the $(Q^{\dagger})_{i,r+2}$ be equal to \een \een \een
\begin{equation}
\frac{ I^{\calX_k}_{+}+R/2}{|\calX^0|+\sum_{j=1}^{r^{\uparrow}-1} I^{\calX_j}_{+}|\calX_{r-j+1}|}(I^{\calX_{r^+}}_{+})^\mathsf{c}.
\end{equation}

\ben \item[]\ben  \item[]\ben \item[iii)] Let all the remaining elements be equal to \een \een \een
\begin{equation} 
\frac{R/2}{|\calX^0|+\sum_{j=1}^{r^{\uparrow}-1} I^{\calX_j}_{+}|\calX_{r-j+1}|}(I^{\calX_{r^+}}_{+})^\mathsf{c}.
\end{equation}

\ben \item[] \ben \item[c)] The elements of all remaining columns are given by 
\ben \item[i)] For all $i\in \calX_k$, $k=1,2,\dots,r^-{-}1$ let 
\een \een \een
\begin{equation}
(Q^{\dagger})_{i,z}=\frac{(I^{\calX_k}_{-})^\mathsf{c}}{|\calX_k|} ,
\end{equation} 
\ben \item[] \ben \item[] \ben \item[]
where $z=1+k$ denotes the $zth$ column. Let all the remaining elements of the $zth$ column be equal to zero. However, if $I^{\calX_k}_{-}=1$, then let all the elements of the $zth$ column be equal with the corresponding elements of the first column, that is,
\een \een \een 
\begin{equation}
(Q^{\dagger})_{1,z}=(Q^{\dagger})_{1,1}, (Q^{\dagger})_{2,z}=(Q^{\dagger})_{2,1},\dots,(Q^{\dagger})_{|\calX|,z}=(Q^{\dagger})_{|\calX|,1}.
\end{equation}
\ben \item[] \ben \item[] \ben \item[ii)]For all $i\in \calX_{r-k+1}$, $k=1,2,\dots,r^+{-}1$ let  \een \een \een
\begin{equation}
(Q^{\dagger})_{i,z}=\frac{( I^{\calX_k}_{+})^\mathsf{c}}{|\calX_k|},
\end{equation} 
\ben \item[] \ben \item[] \ben \item[]
where $z=r+2-k$ denotes the $zth$ column. Let all the remaining elements of the $zth$ column be equal to zero. However, if $ I^{\calX_k}_{+}=1$, then let all the elements of the $zth$ column be equal with the corresponding elements of the last column, that is,
\een \een \een 
\begin{equation}
(Q^{\dagger})_{1,z}=(Q^{\dagger})_{1,|\calX|}, (Q^{\dagger})_{2,z}=(Q^{\dagger})_{2,|\calX|},\dots,(Q^{\dagger})_{|\calX|,z}=(Q^{\dagger})_{|\calX|,|\calX|}.
\end{equation}
\end{algorithm}

Once the $Q^\dagger$ matrix is constructed, as a function of the TV parameter $R$, then by \eqref{MC2} the solution of optimization \eqref{max.entropy} is readily available, and hence, by Definition \ref{lower.dim.process2}, the lower dimensional process $\{Y_t:t=0,1\dots\}$ with invariant distribution $\bar{\nu}$  is obtained, either by adding all equal elements of $\nu^*\in \mathbb{P}(\calX)$, or by defining a $Q$ matrix to be equal to $Q^{\dagger}$, after the merging of all equal columns (by adding them). Hence \begin{equation}
\bar{\nu}=\mu Q=\mu PQ ,
\end{equation}
where the dimensions of $Q$ matrix are based on the value of the TV parameter $R\in[0,2]$.

Before we proceed with the solution of \eqref{KLentropy}, we provide a simple example in order to explain each step of Algorithm \ref{algorithm2}.

\begin{example}\label{ex.alg.entropy}
Let $\mu=[\mu_1 \ \mu_2\  \mu_3\ \mu_4]$, where $\mu_1>\mu_2>\mu_3>\mu_4$, and also assume that $\mu_1>\mu_2>\frac{1}{|\calX|}$ and $\mu_4<\mu_3<\frac{1}{|\calX|}$, where $|\calX|=4$. For simplicity of presentation it is assumed that the optimum probabilities $\nu_i^*$, $i\in \calX$, as a function of $R$ are as shown in Fig.\ref{fig100}. 

Initialization step. For $R=0$, and from Remark \ref{remalg2}, we conclude that $\ell_1<\ell_2<\ell_3<\ell_4$, and therefore the support sets are equal to $\calX^0=\{4\}$, $\calX_0=\{1\}$, $\calX_1=\{2\}$ and $\calX_2=\{3\}$. The number of the $\calX_k$ sets is equal to $r=2$. The number of $\mu_i$, $i\in \calX$, which are greater (or equal) than $\frac{1}{|\calX|}=0.25$ (and also $\mu_i \neq \mu_j$, $i,j\in \calX$) is $r^-=2$. Similarly, the number of $\mu_i$ which are strictly smaller than $\frac{1}{|\calX|}=0.25$ (and also not equal to each other) is also $r^+=2$.

Step.1 From \eqref{alg2ind1}-\eqref{alg2ind2}, the indicator functions $I_{-}^{\calX_1}$ and $I_{-}^{\calX_2}$ are given by
\begin{equation*}
I_{-}^{\calX_1}{\triangleq} \left\{ \,
\begin{IEEEeqnarraybox}[][c]{l?s}
\IEEEstrut
1, & if $\mu_1{-}\frac{R}{2}{\leq} \mu_2$,\\
0, & otherwise,
\IEEEstrut
\end{IEEEeqnarraybox}
\right.
\qquad
I_{-}^{\calX_2}{\triangleq} \left\{ \,
\begin{IEEEeqnarraybox}[][c]{l?s}
\IEEEstrut
1, & if $\frac{\mu_1+\mu_2-R/2}{2}{\leq} 0.25$,\\
0, & otherwise,
\IEEEstrut
\end{IEEEeqnarraybox}
\right.
\end{equation*}
and from \eqref{alg2ind3}-\eqref{alg2ind4}, the indicator functions $I_{+}^{\calX_1}$ and $I_{+}^{\calX_2}$ are given by
\begin{equation*}
I_{+}^{\calX_1}{\triangleq} \left\{ \,
\begin{IEEEeqnarraybox}[][c]{l?s}
\IEEEstrut
1, & if $\mu_4{+}\frac{R}{2}{\geq} \mu_3$,\\
0, & otherwise,
\IEEEstrut
\end{IEEEeqnarraybox}
\right.
\qquad
I_{+}^{\calX_2}{\triangleq} \left\{ \,
\begin{IEEEeqnarraybox}[][c]{l?s}
\IEEEstrut
1, & if $\frac{\mu_3+\mu_4+R/2}{2}{\geq} 0.25$,\\
0, & otherwise.
\IEEEstrut
\end{IEEEeqnarraybox}
\right.
\end{equation*}
The values of the indicator functions for $R\in [0,2]$ are shown in Fig.\ref{fig100}. For $0 \leq R< R_1$, that is, before a merge occurs, all indicator functions are equal to zero. If a merge occurs the respective indicator functions become equal to one, until for some $R\geq R_3$, where all indicator functions are equal to one.

Step.2 Let $Q^{\dagger}$ be an $4\times 4$ matrix. For $0 \leq R< R_1$,
\bes Q^{\dagger}=\left(
         \begin{array}{cccc}
           1-R/2 & 0 & 0 & R/2 \\
           -R/2 & 1 & 0 & R/2 \\
          -R/2 & 0 & 1 &R/2 \\
           -R/2 & 0 & 0 &1+R/2 \\
         \end{array}
       \right),\ees
and since no equal columns exist then $Q^{\dagger}=Q$. For $R_1\leq R< R_2$,
\bes Q^{\dagger}=\left(
         \begin{array}{cccc}
           \frac{1-R/2}{2} & \frac{1-R/2}{2} & 0 & R/2 \\
          \frac{1 -R/2}{2} & \frac{1-R/2}{2} & 0 & R/2 \\
          -R/4 & -R/4 & 1 &R/2 \\
           -R/4 & -R/4 & 0 &1+R/2\\
         \end{array}
       \right)\Longrightarrow Q=\left(
         \begin{array}{ccc}
           1-R/2 & 0& R/2 \\
           1-R/2 &0& R/2  \\
           -R/2 &1& R/2  \\
           -R/2 &0& 1+R/2 \\
         \end{array}
       \right).\ees         
       
       \noi For $R_2\leq R< R_3$,
\bes Q^{\dagger}=\left(
         \begin{array}{cccc}
           \frac{1-R/2}{2} & \frac{1-R/2}{2} & R/4 & R/4 \\
          \frac{1 -R/2}{2} & \frac{1-R/2}{2} & R/4& R/4 \\
          -R/4 & -R/4 & \frac{1+R/2}{2} &\frac{1+R/2}{2} \\
           -R/4 & -R/4 & \frac{1+R/2}{2} &\frac{1+R/2}{2}\\
         \end{array}
       \right)\Longrightarrow Q=\left(
         \begin{array}{cc}
           1-R/2& R/2 \\
           1-R/2 & R/2  \\
           -R/2 & 1+R/2  \\
           -R/2 & 1+R/2  \\
         \end{array}
       \right).\ees

For $R\geq R_3$,
\begin{equation*} Q^{\dagger}=\left(
         \begin{array}{cccc}
           0.25 & 0.25 & 0.25 & 0.25 \\
           0.25 & 0.25 & 0.25 & 0.25 \\
          0.25 & 0.25 & 0.25 &0.25 \\
          0.25 & 0.25 & 0.25 &0.25 \\
         \end{array}
       \right)\Longrightarrow Q=\left(
         \begin{array}{cccc}
          1  \\
           1   \\
          1 \\
           1  \\
         \end{array}
       \right).
\end{equation*}

\begin{figure}[htbp!]
\centering
\includegraphics[width=0.95\columnwidth]{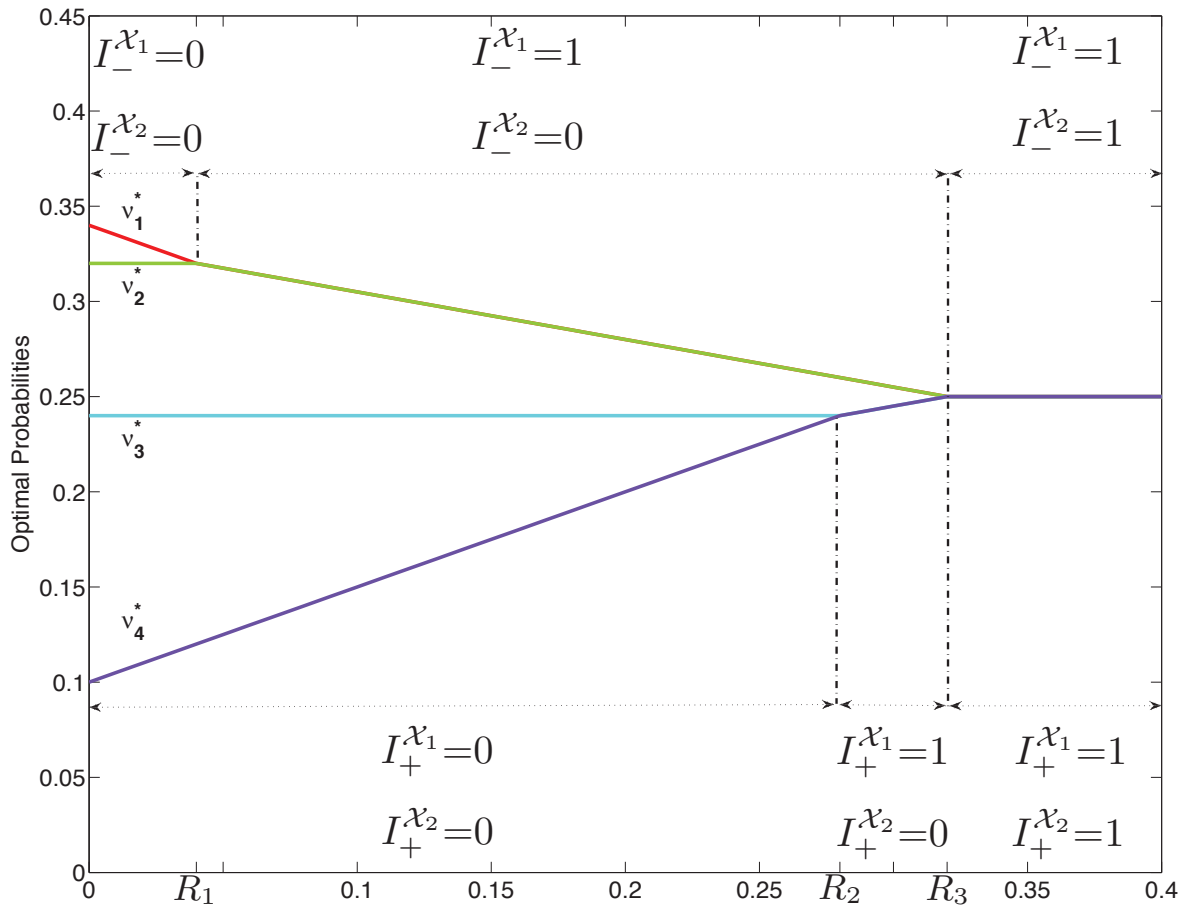}
\caption[]{Optimal Probabilities as a function of $R$.}
\label{fig100}
\end{figure}
Note that, the dimension of matrix $Q$ is based on the value of total variation distance parameter $R$. For $0< R\leq R_1$ its dimension is equal to $(|\calX|)\times (2+r)$. Whenever two columns become equal (that is, an indicator function is activated) they are merged, until for some $R\geq R_2$, where matrix $Q$ is transformed into column vector of dimension $(|\calX|)\times (1)$. 
\end{example}

 Next, we proceed with the solution of \eqref{KLentropy}, by letting $\bar{\nu}$ to denote the invariant distribution of a lower dimensional Markov process $(\bar{\nu},\Phi)$. To this end, we next define an optimal partition function for the approximation problem, based on maximum entropy principle at values of TV parameter $R$, for which an aggregation of the states occurs (i.e., see Example \ref{ex.alg.entropy}, Fig.\ref{fig100}, for values of $R=R_1$, $R_2$ and $R_3$.).

\begin{definition}\label{part.function.entr}(partition function) Let $\calX$ and $\calY$ be two finite dimensional state-spaces with $|\calY|<|\calX|$. Define a surjective (partition) function $\varphi:\calX\longmapsto \calY$ as follows \begin{equation}
\forall i,j\in\calX,\quad \varphi(i)=\varphi(j)=k\in\calY\quad \mbox{if}\quad \nu^*_i=\nu^*_j.
\end{equation}
\end{definition}

Note that, once the optimal probabilities $\nu_i^*$, $\forall i\in\calX$ are obtained, we can easily identify the values of $R$ for which an aggregation of the states occurs. Next, we reproduce the main theorem of \cite{Meyn09}, which gives the solution of $\Phi$ that solves \eqref{KLentropy}.
\begin{theorem}
Let $(\mu,P,\calX)$ be a FSM process and $\varphi$ be the partition function of Definition \ref{part.function.entr}. For optimization \eqref{KLentropy}, the solution of $\Phi$ is given by
\begin{equation}
\Phi_{kl}=\frac{u^{(k)}\Pi Pu^{(\ell)'}}{\bar{\nu}_k},\quad k,\ell\in \calY
\end{equation}
where $\Pi=diag(\nu^*)$, $u^{(k)'}$ is the transpose of $u^{(k)}$, and $u^{(k)}$ is a $1\times|\calX|$ row vector defined by
\begin{equation} u_i^{(k)} = \left\{
  \begin{array}{l l}
    1 & \quad \text{if $\varphi(i)=k$}\\
   0 & \quad \text{otherwise}
  \end{array} \right.\end{equation}
\end{theorem}

\begin{IEEEproof}
See \cite{Deng:11}.
\end{IEEEproof}

\section{Examples}\label{sec.examples}

\subsection{Markov chain approximation with a small number of states}\label{example1}
In this example, we employ the theoretical results obtained in preceding sections to approximate a $4$-state FSM process $(\mu,P,\calX)$ with transition probability matrix given by \bea P=\left[
         \begin{array}{cccc}
           0.4 & 0.2& 0.3  & 0.1 \\
           0.3 & 0.5& 0.1  & 0.1 \\
           0.2 & 0.3& 0.4  &0.1 \\
           0.6 & 0.2& 0.1  & 0.1 \\
         \end{array}
       \right] ,\eea
and steady state nominal probability vector   equal to \bea\mu=[0.34 \hso 0.32 \hso 0.24\hso  0.1].\eea
In particular, in Section \ref{ex1.3}, we solve approximation problem based on \emph{Method 1}. In Section \ref{ex1.1} we solve the approximation problem based on occupancy distribution, and in Section \ref{ex1.2} based on entropy principle of \emph{Method 2}.

\subsubsection{Solution of Problem \ref{problem3}}\label{ex1.3}
Let $\ell=\{\ell\in\mathbb{R}_+^4:\ell_1>\ell_2>\ell_3>\ell_4\}$, then the support sets are given by $\calX^0=\{1\}$, $\calX_0=\{4\}$, $\calX_1=\{3\}$ and $\calX_2=\{2\}$, and by \eqref{eq.55d}, $R_{\max,1}=1.2$, $R_{\max,2}=1.4$, $R_{\max,3}=1.6$ and $R_{\max,4}=0.8$. By employing Theorem \ref{MCthmocc}, the optimal $\Phi^\dagger$ and $\Phi$ matrices are obtained as a function of TV parameter $R$, as shown in Table \ref{Table.Ex.c}. Note that, in contrast with Problems \ref{problem.occ}-\ref{problem.entr}, where the approximation is performed only for values of $R$ for which a reduction of the states occurs, the solution of Problem \ref{problem3} is obtained for all values of total variation parameter.

\begin{table}[!h]\setlength{\tabcolsep}{2pt}\centering
\begin{tabular}{|c ||c|| c|  }
\hline
 {\textsc{$R$}} & \textsc{$\Phi^{\dagger}$} & \textsc{$\Phi$}\\
 \hline \hline
0& $\left[\footnotesize {\begin{tabular}{cccc}
   .4 & .2 & .3 & .1 \\
   .3 & .5 & .1 & .1 \\
   .2 & .3 & .4 & .1 \\
   .6 & .2 & .1 & .1 \\
 \end{tabular}}\right]$
 &$\left[\footnotesize {\begin{tabular}{cccc}
   .4 & .2 & .3 & .1 \\
   .3 & .5 & .1 & .1 \\
   .2 & .3 & .4 & .1 \\
   .6 & .2 & .1 & .1 \\
 \end{tabular}}\right]$\\ \hline \hline
 0.2& $\left[\footnotesize{\begin{tabular}{cccc}
   .5 & .2 & .3 &0\\
   .4 & .5 & .1  &0\\
   .3 & .3 & .4  &0\\
   .7 & .2 & .1  &0\\
 \end{tabular}}\right]$
&$\left[\footnotesize{\begin{tabular}{ccc}
   .5 & .2 & .3  \\
   .4 & .5 & .1  \\
   .3 & .3 & .4  \\
 \end{tabular}}\right]$\\\hline \hline
 1& $\left[\footnotesize{\begin{tabular}{cccc}
  .9 & .1 &0 &0 \\
   .8 & .2 &0 &0\\
   .7 & .3 &0 &0\\
   1 & 0 &0 &0\\
 \end{tabular}}\right]$
&$\left[\footnotesize{\begin{tabular}{cccc}
   .9 & .1  \\
   .8 & .2 \\
 \end{tabular}}\right]$\\\hline \hline
 1.4& $\left[\footnotesize{\begin{tabular}{cccc}
   1 &0 &0 &0 \\
   1 &0 &0 &0 \\
   1 &0 &0 &0 \\
   1 &0 &0 &0 \\
 \end{tabular}}\right]$
&$\left[\footnotesize{\begin{tabular}{cccc}
   1  \\
 \end{tabular}}\right]$\\ \hline
\end{tabular}
\caption{\small{Optimal results obtained by the  Approximation based on \emph{Method 1}.}}
\label{Table.Ex.c}
\end{table}

\subsubsection{Solution of Problem \ref{problem.occ}}\label{ex1.1}
By employing Algorithm \ref{algorithm1}, with $\ell_i\tri\mu_i$, $i=1,\dots,4$, and support sets given by $\calX^0=\{1\}$, $\calX_0=\{4\}$, $\calX_1=\{3\}$ and $\calX_2=\{2\}$ the maximizing distribution of \eqref{max.occupancy} exhibits a water-filling behavior as depicted in Fig.\ref{fig1000}. For values of TV parameter $0\leq R\leq R_1=0.2$, all maximizing probabilities $\nu_i^*$, $i=1,\dots,4$, are greater than zero and hence $|\calY|=4=|\calX|$ and $\bar{\nu}_i=\nu^*_i$, $i=1,\dots,4$. However, for $R_1\leq R<R_2=0.68$, $|\calY|=3<|\calX|=4$ since $\nu^*_4$ becomes equal to zero and hence $\bar{\nu_i}=\nu_i^*$, $i=1,2,3$. The procedure follows until for some $R\geq R_3=1.32$ in which $|\calY|=1$ and $\bar{\nu}_1=\nu_1^*=1$.

From the above discussion, it is clear that, the solution of approximation problem based on occupancy distribution is described via a water-filling deletion of states with the smallest invariant probability and maintaining and strengthening the states with the highest invariant probability, and hence a lower dimensional distribution $\bar{\nu}$ is obtained which is then applied to the problem of Markov by Markov approximation. For the solution of \eqref{KLoccupancy}, first we find an optimal partition function $\varphi$ and then we calculate a transition probability matrix $\Phi$ which best approximates transition matrix $P$ only for values of $R$ for which a reduction of states occurs, that is, for $R=0,~0.2,~0.68$ and $1.32$. The optimal results are depicted in Table \ref{Table.Ex.a}.

\begin{table}[!h]\setlength{\tabcolsep}{2pt}\centering
\begin{tabular}{|c|| c ||c|| c| | c| }
\hline
 {\textsc{$R$}}& {\textsc{$\bar{\nu}$}} & \textsc{$Q$} & \textsc{$\varphi$} & \textsc{$\Phi$}\\
\hline \hline
0& \footnotesize{\begin{tabular}{cccc}
    [.34 & .32 & .24 & .1]\\
  \end{tabular}}
&$\left[\footnotesize {\begin{tabular}{cccc}
   1 & 0 & 0 & 0 \\
   0 & 1 & 0 & 0 \\
  0 & 0 & 1 & 0 \\
   0 & 0 & 0 & 1 \\
 \end{tabular}}\right]$
 &\begin{tabular}{c}
 $\varphi(1)=1$ \\  
 $\varphi(2)=2$  \\ 
 $\varphi(3)=3$  \\ 
 $\varphi(4)=4$  \\ 
 \end{tabular}
 &$\left[\footnotesize {\begin{tabular}{cccc}
   .4 & .2 & .3 & .1 \\
   .3 & .5 & .1 & .1 \\
  .2 & .3 & .4 & .1 \\
   .6 & .2 & .1 & .1 \\
 \end{tabular}}\right]$\\
\hline \hline
 0.2& \footnotesize{\begin{tabular}{cccc}
    [.44 & .32 & .24]\\
  \end{tabular}}
  &$\left[\footnotesize{\begin{tabular}{ccc}
   1.1 & 0 & -.1  \\
   .1 & 1 & -.1  \\
   .1 & 0 & .9  \\
   .1 & 0 & .9  \\
 \end{tabular}}\right]$
&\begin{tabular}{c}
 $\varphi(1)=1$ \\  
 $\varphi(2)=2$  \\ 
 $\varphi(3)=3$  \\ 
 $\varphi(4)=1$  \\ 
 \end{tabular}
&$\left[\footnotesize {\begin{tabular}{ccc}
   .5455 & .2 & .2545  \\
   .4 & .5 & .1  \\
  .3 & .3 & .4  \\
 \end{tabular}}\right]$\\ 
 \hline \hline
 0.68& \footnotesize{\begin{tabular}{cccc}
    [0.68 & 0.32]\\
  \end{tabular}}
  &$\left[\footnotesize{\begin{tabular}{cccc}
  1.34 & -.34  \\
   .34 & .66 \\
   .34 & .66 \\
   .34 & .66 \\
 \end{tabular}}\right]$
&\begin{tabular}{c}
 $\varphi(1)=1$ \\  
 $\varphi(2)=2$  \\ 
 $\varphi(3)=1$  \\ 
 $\varphi(4)=1$  \\ 
 \end{tabular}
&$\left[\footnotesize {\begin{tabular}{cc}
   .7647 & .2353  \\
   .5 & .5  \\
 \end{tabular}}\right]$\\ 
 \hline \hline
 1.32& \footnotesize{\begin{tabular}{cccc}
    [1]\\
  \end{tabular}}
  &$\left[\footnotesize{\begin{tabular}{cccc}
   2.94  \\
   0  \\
   0  \\
   0  \\
 \end{tabular}}\right]$
&\begin{tabular}{c}
 $\varphi(1)=1$ \\  
 $\varphi(2)=1$  \\ 
 $\varphi(3)=1$  \\ 
 $\varphi(4)=1$  \\ 
 \end{tabular}
&$\left[\footnotesize {\begin{tabular}{c}
   1 
 \end{tabular}}\right]$\\ 
  \hline
\end{tabular}
\caption{\small{Optimal results obtained by the  Approximation based on occupancy distribution\vspace{-0.3cm}.}}
\label{Table.Ex.a}
\end{table}

\subsubsection{Solution of Problem \ref{problem.entr}}\label{ex1.2}
By employing Algorithm \ref{algorithm2}, with $\ell_i\tri-\log\nu_i$, $i=1,\dots,4$, the support sets are calculated for $R=0$, where $\nu_i^*=\mu_i$ and hence $\ell_i=-\log\mu_i$, and are equal to $\calX^0=\{4\}$, $\calX_0=\{1\}$, $\calX_1=\{2\}$ and $\calX_2=\{3\}$. The maximizing distribution of \eqref{max.entropy} exhibits a water-filling like behavior as depicted in Fig.\ref{fig100}. For values of $0\leq R<R_1=0.04$, $|\calY|=4=|\calX|$ since $\nu_i^*\neq \nu_j^*$ for $i\neq j$, $i,j=1,\dots,4$ and hence $\bar{\nu}_i=\nu_i^*$, $i=1,\dots,4$. For $R_1\leq R<R_2=0.28$, $|\calY|=3<|\calX|=4$ since $\nu_1^*$ becomes equal to $\nu_2^*$ and hence $\bar{\nu}_1=\nu_1^*+\nu_2^*$ and $\bar{\nu}_i=\nu_i^*$, $i=3,4$. The procedure follows until for some $R\geq R_3=0.32$ in which $|\calY|=1$ and $\bar{\nu}_1=\sum_{i=1}^4\nu_i^*=\frac{1}{4}$.

In summary, the solution of approximation problem based on entropy principle is described via aggregation of states, that is, by grouping certain states of the original Markov chain to obtain the approximating reduced state process.Then the lower dimensional distribution $\bar{\nu}$ is applied to problem \eqref{KLentropy}. The optimal partition function $\varphi$ and the transition probability matrix $\Phi$ which minimizes the KL divergence rate for values of $R=0,~0.04,~0.28$ and $0.32$ are as shown in Table \ref{Table.Ex.b}.
\begin{figure*}[!htbp]
\centering
\subfloat[][]{
\label{fig1.1} 
\includegraphics[ width=0.32\linewidth]{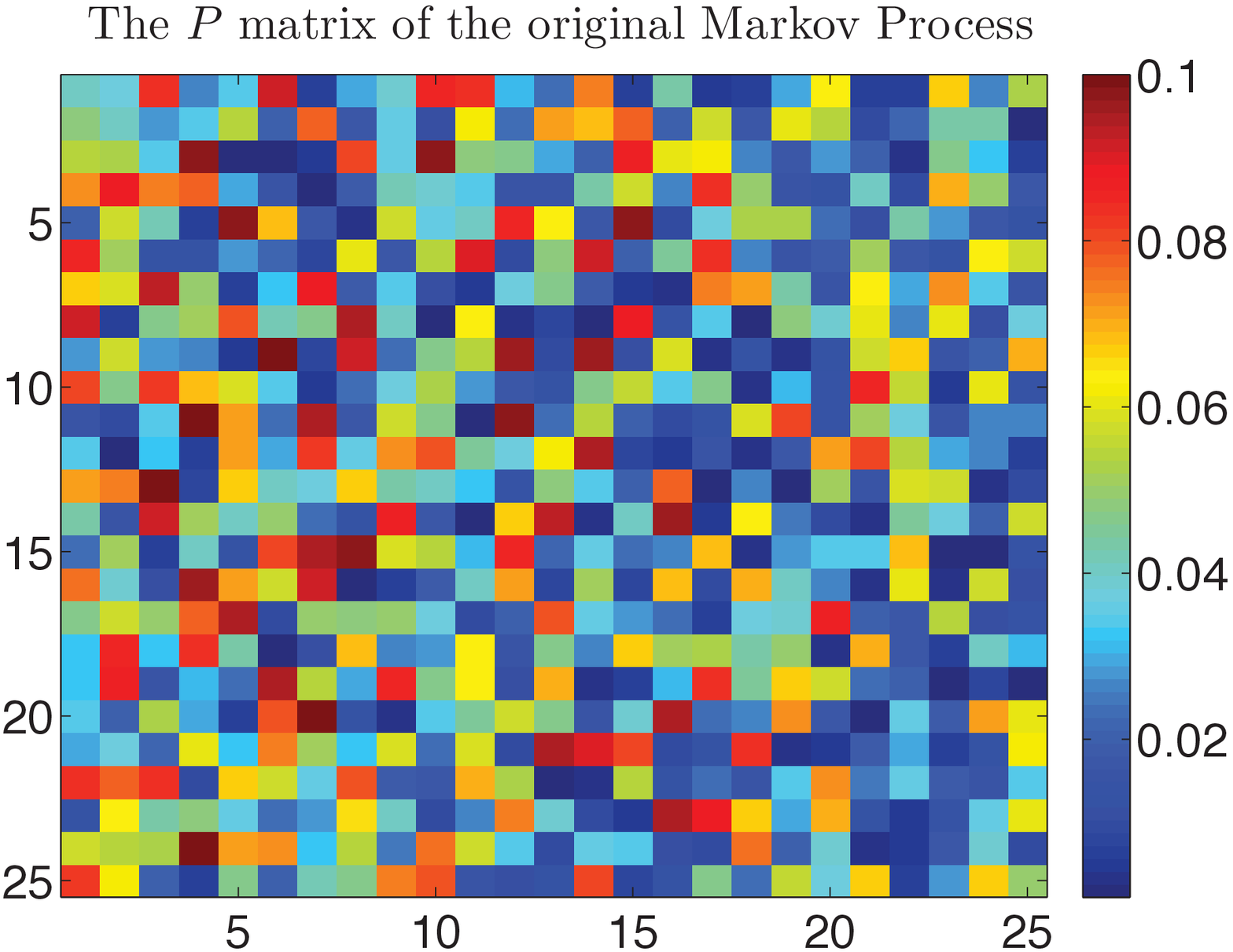}}
\subfloat[][]{
\label{fig1.2} 
\includegraphics[width=0.32\linewidth]{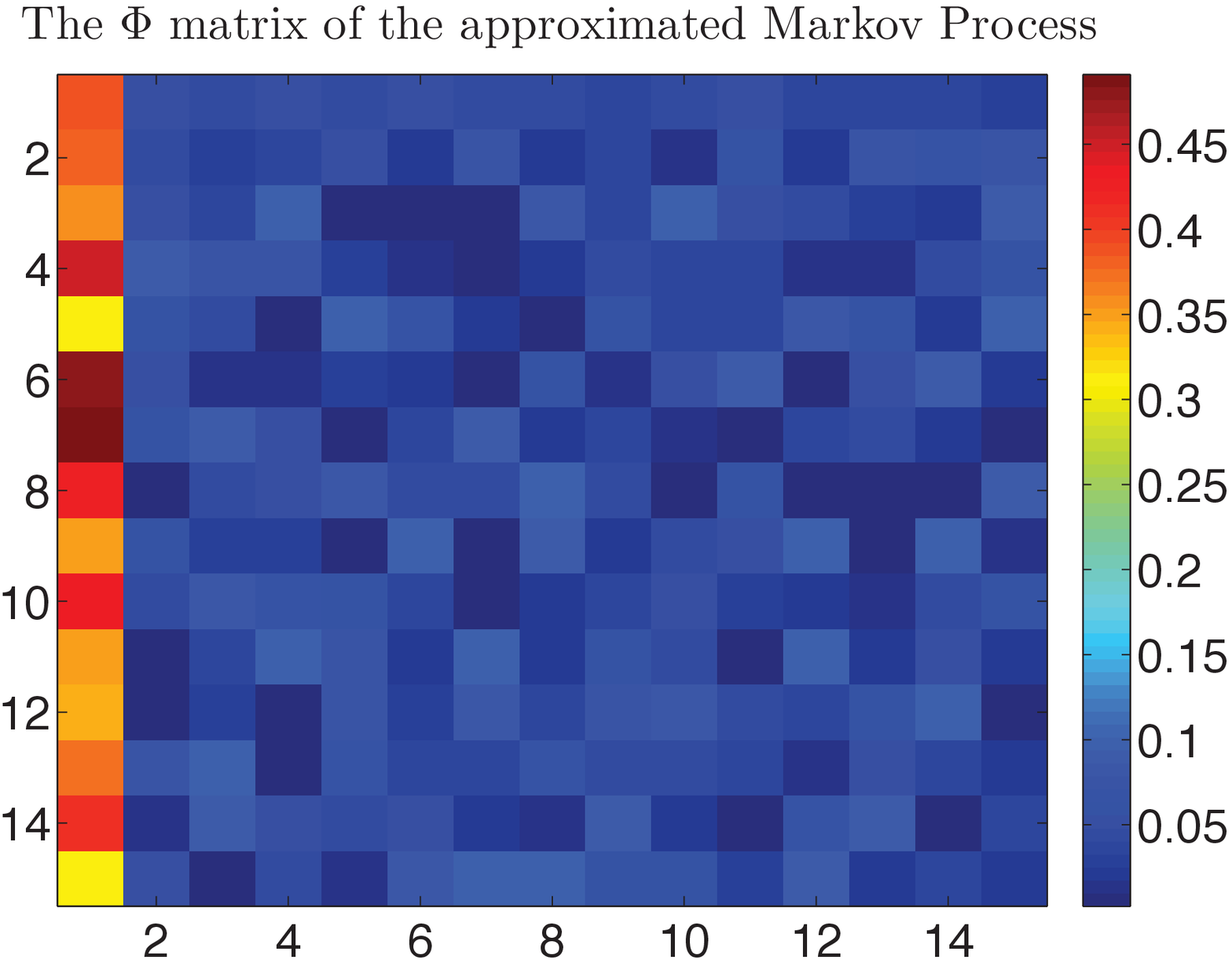}}
\subfloat[][]{
\label{fig1.3} 
\includegraphics[ width=0.32\linewidth]{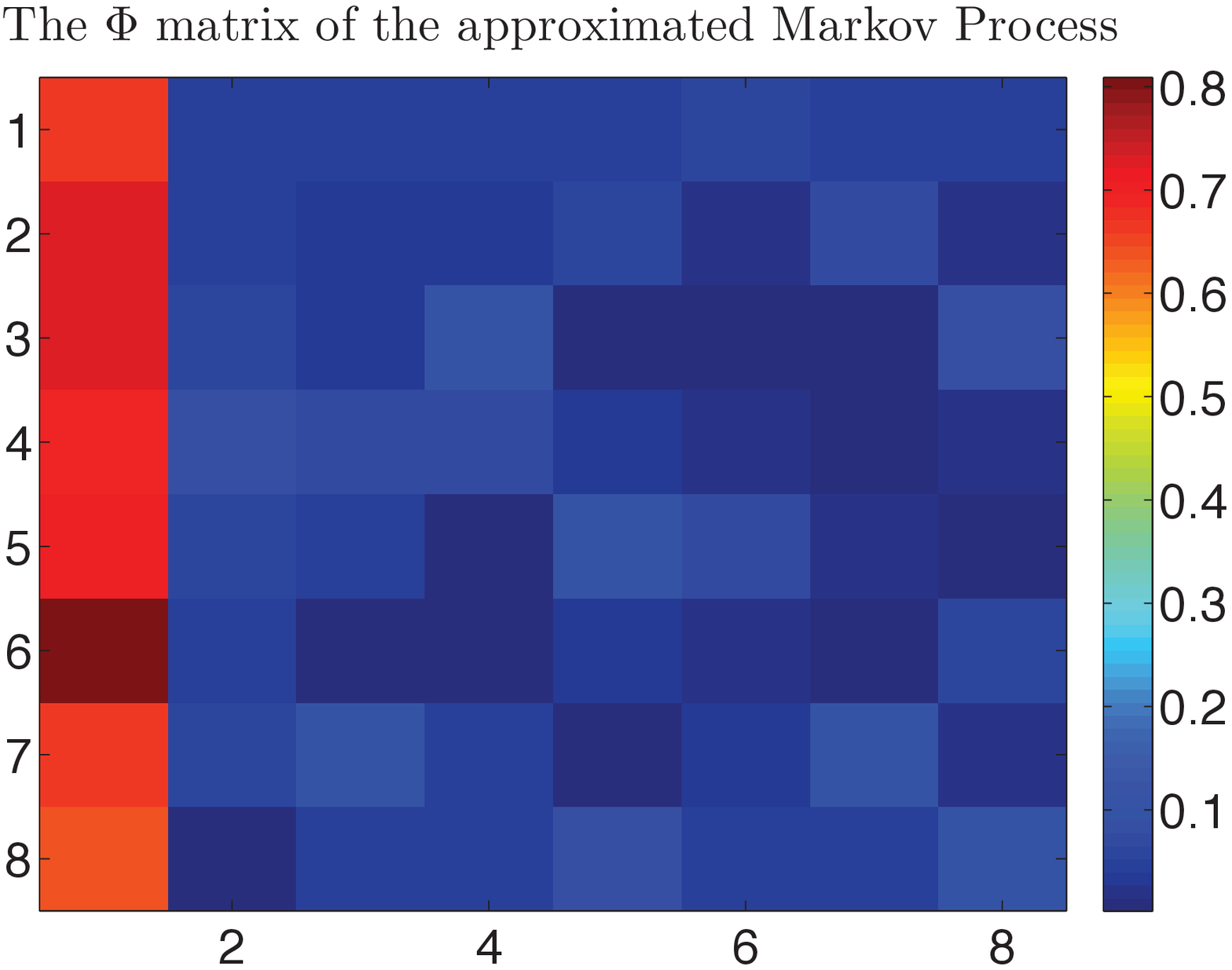}}
\hfil
\vspace{-.3cm}
\subfloat[][]{
\label{fig1.4} 
\includegraphics[ width=0.32\linewidth]{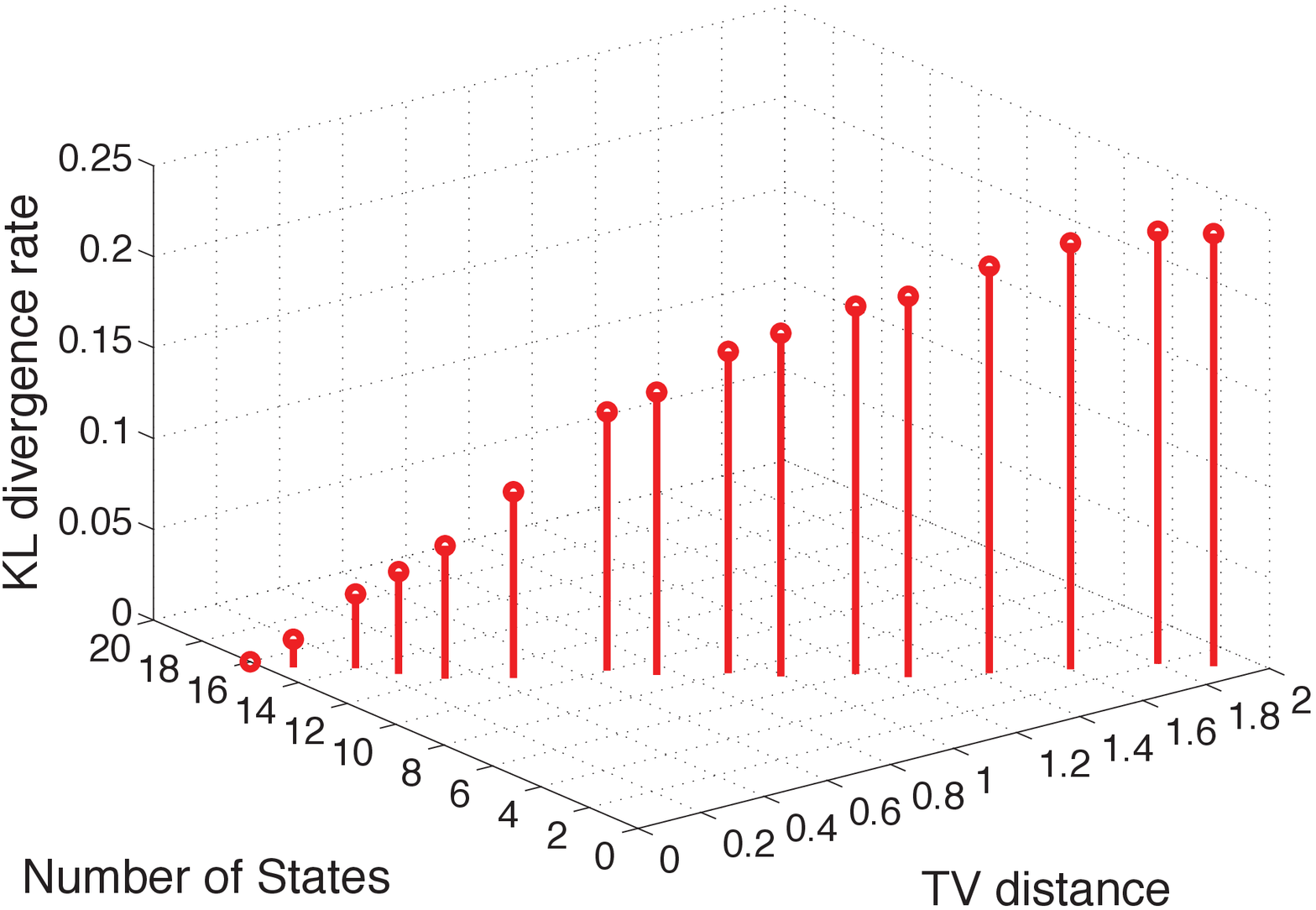}}
\subfloat[][]{
\label{fig1.5} 
\includegraphics[ width=0.32\linewidth]{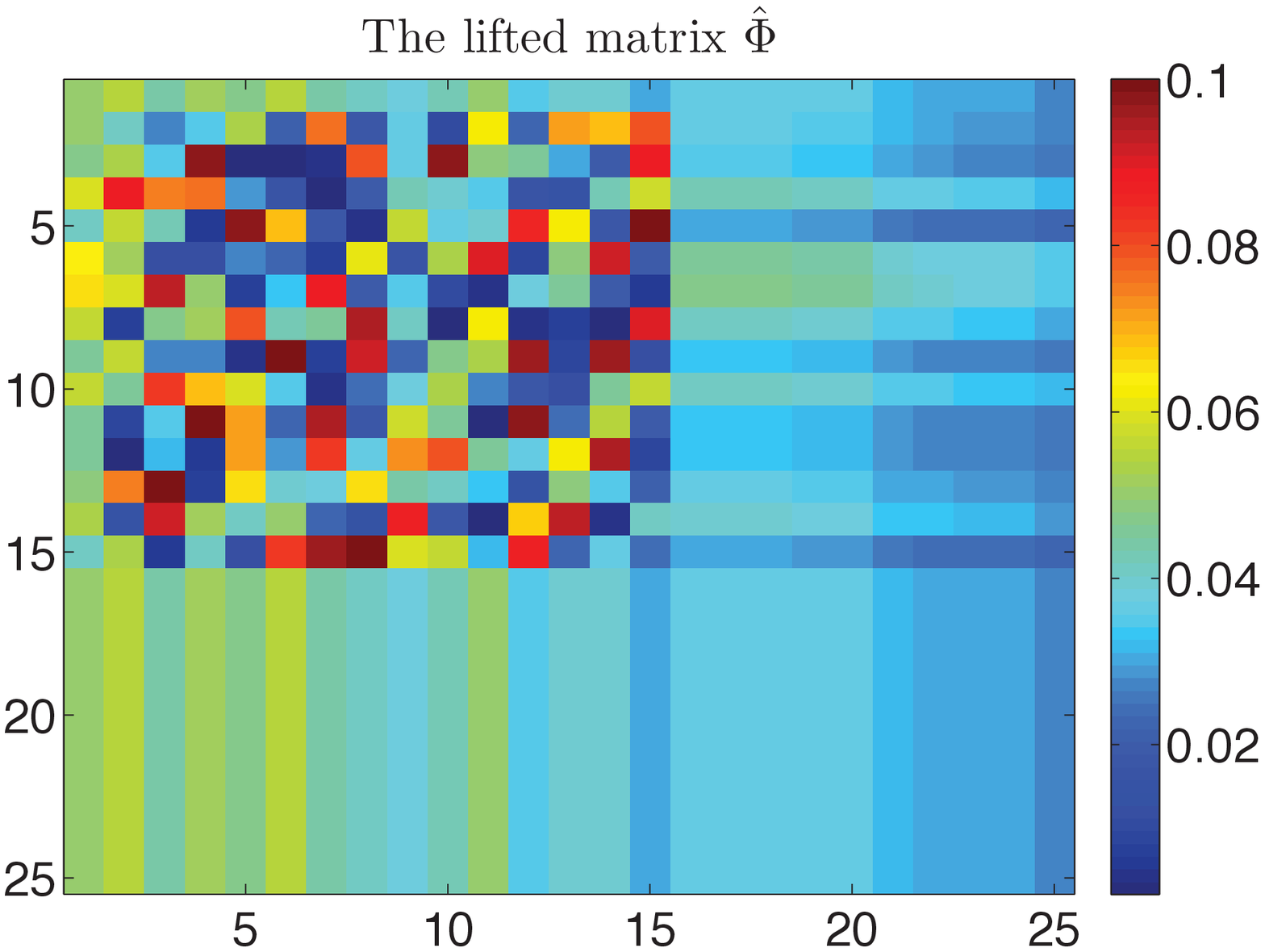}}
\subfloat[][]{
\label{fig1.6} 
\includegraphics[ width=0.32\linewidth]{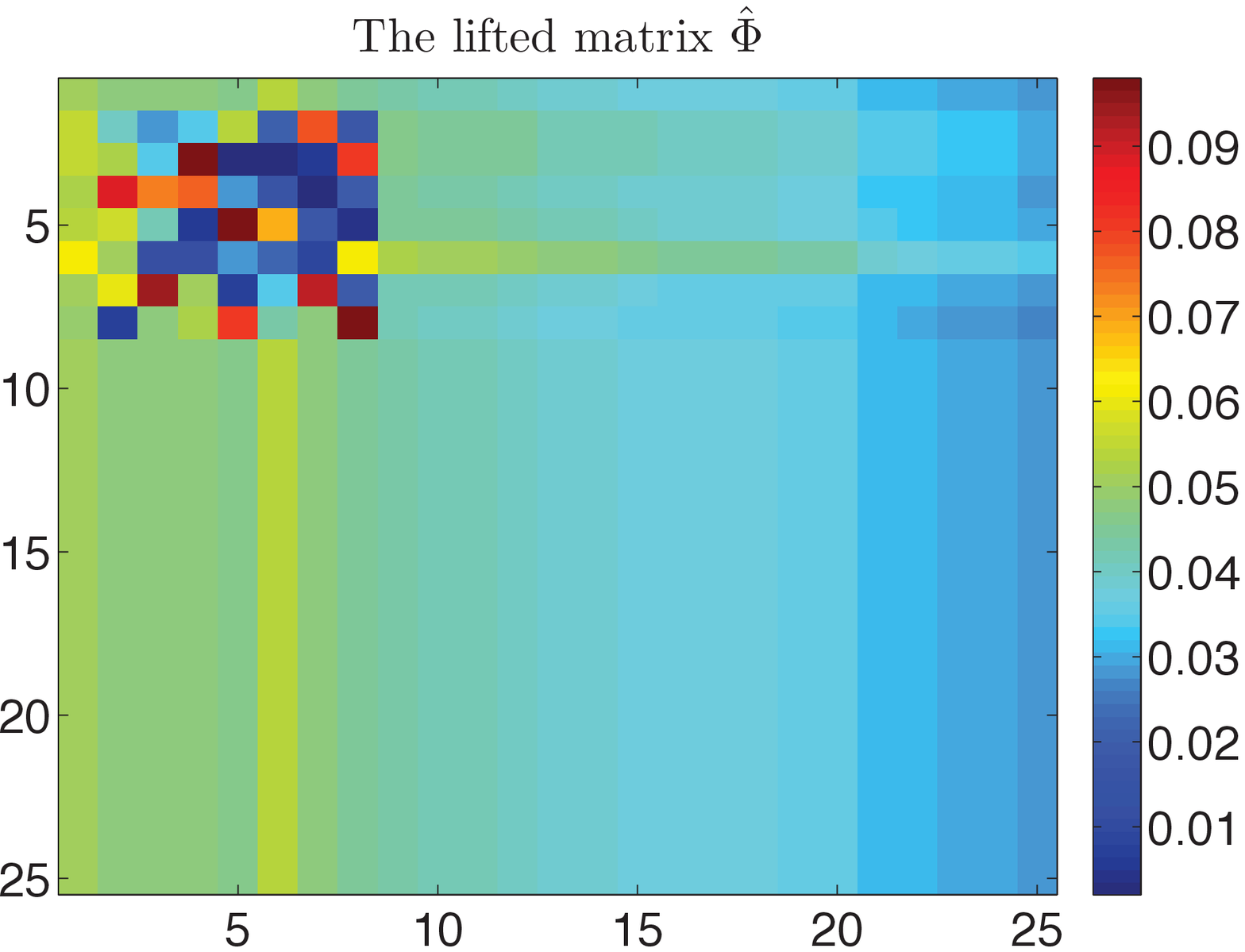}}
\caption[Optimal Solution of Example B]{Approximation results based on occupancy distribution: Plot (a) depicts the $P$ matrix of the original Markov process. Plot (b) depicts a $15$-state approximation. Plot (c) depicts an $8$-state approximation. Plot (d) depicts the KL divergence rate. Plot (e) depicts the lifted $\hat{\Phi}$ matrix for the $15$-state approximation. Plot (f) depicts the lifted $\hat{\Phi}$ matrix for the $8$-state approximation.}
\label{fig1}
\end{figure*}

\begin{table}[!h]\setlength{\tabcolsep}{2pt}\centering
\begin{tabular}{|c|| c ||c|| c| | c| }
\hline
 {\textsc{$R$}}& {\textsc{$\bar{\nu}$}} & \textsc{$Q$} & \textsc{$\varphi$} & \textsc{$\Phi$}\\
\hline \hline
0& \footnotesize{\begin{tabular}{cccc}
    [.34 & .32 & .24 & .1]\\
  \end{tabular}}
&$\left[\footnotesize {\begin{tabular}{cccc}
   1 & 0 & 0 & 0 \\
   0 & 1 & 0 & 0 \\
  0 & 0 & 1 & 0 \\
   0 & 0 & 0 & 1 \\
 \end{tabular}}\right]$
 &\begin{tabular}{c}
 $\varphi(1)=1$ \\  
 $\varphi(2)=2$  \\ 
 $\varphi(3)=3$  \\ 
 $\varphi(4)=4$  \\ 
 \end{tabular}
 &$\left[\footnotesize {\begin{tabular}{cccc}
   .4 & .2 & .3 & .1 \\
   .3 & .5 & .1 & .1 \\
  .2 & .3 & .4 & .1 \\
   .6 & .2 & .1 & .1 \\
 \end{tabular}}\right]$\\
\hline \hline
 0.04& \footnotesize{\begin{tabular}{cccc}
    [.64 & .24 & .12]\\
  \end{tabular}}
  &$\left[\footnotesize{\begin{tabular}{ccc}
   .98 & 0 & .02  \\
   .98 & 0 & .02  \\
   -.02 & 1 & .02  \\
   -.02 & 0 & 1.02  \\
 \end{tabular}}\right]$
&\begin{tabular}{c}
 $\varphi(1)=1$ \\  
 $\varphi(2)=1$  \\ 
 $\varphi(3)=2$  \\ 
 $\varphi(4)=3$  \\ 
 \end{tabular}
&$\left[\footnotesize {\begin{tabular}{ccc}
   .7 & .2 & .1  \\
   .5 & .4 & .1  \\
  .8 & .1 & .1  \\
 \end{tabular}}\right]$\\ 
 \hline \hline
 0.28& \footnotesize{\begin{tabular}{cccc}
    [0.52 & 0.48]\\
  \end{tabular}}
  &$\left[\footnotesize{\begin{tabular}{cccc}
  .86 & .14  \\
   .86 & .14 \\
   -.14 & 1.14 \\
   -.14 & 1.14 \\
 \end{tabular}}\right]$
&\begin{tabular}{c}
 $\varphi(1)=1$ \\  
 $\varphi(2)=1$  \\ 
 $\varphi(3)=2$  \\ 
 $\varphi(4)=2$  \\ 
 \end{tabular}
&$\left[\footnotesize {\begin{tabular}{cc}
   .7 & .3  \\
   .65 & .35  \\
 \end{tabular}}\right]$\\ 
 \hline \hline
 0.32& \footnotesize{\begin{tabular}{cccc}
    [1]\\
  \end{tabular}}
  &$\left[\footnotesize{\begin{tabular}{cccc}
   1  \\
   1  \\
   1  \\
   1  \\
 \end{tabular}}\right]$
&\begin{tabular}{c}
 $\varphi(1)=1$ \\  
 $\varphi(2)=1$  \\ 
 $\varphi(3)=1$  \\ 
 $\varphi(4)=1$  \\ 
 \end{tabular}
&$\left[\footnotesize {\begin{tabular}{c}
   1 
 \end{tabular}}\right]$\\ 
  \hline
\end{tabular}
\caption{\small{Optimal results obtained by the  Approximation based on entropy principle\vspace{-0.3cm}.}}
\label{Table.Ex.b}
\end{table}

\begin{figure*}[t]
\centering
\subfloat[][]{
\label{fig2.1} 
\includegraphics[ width=0.32\linewidth]{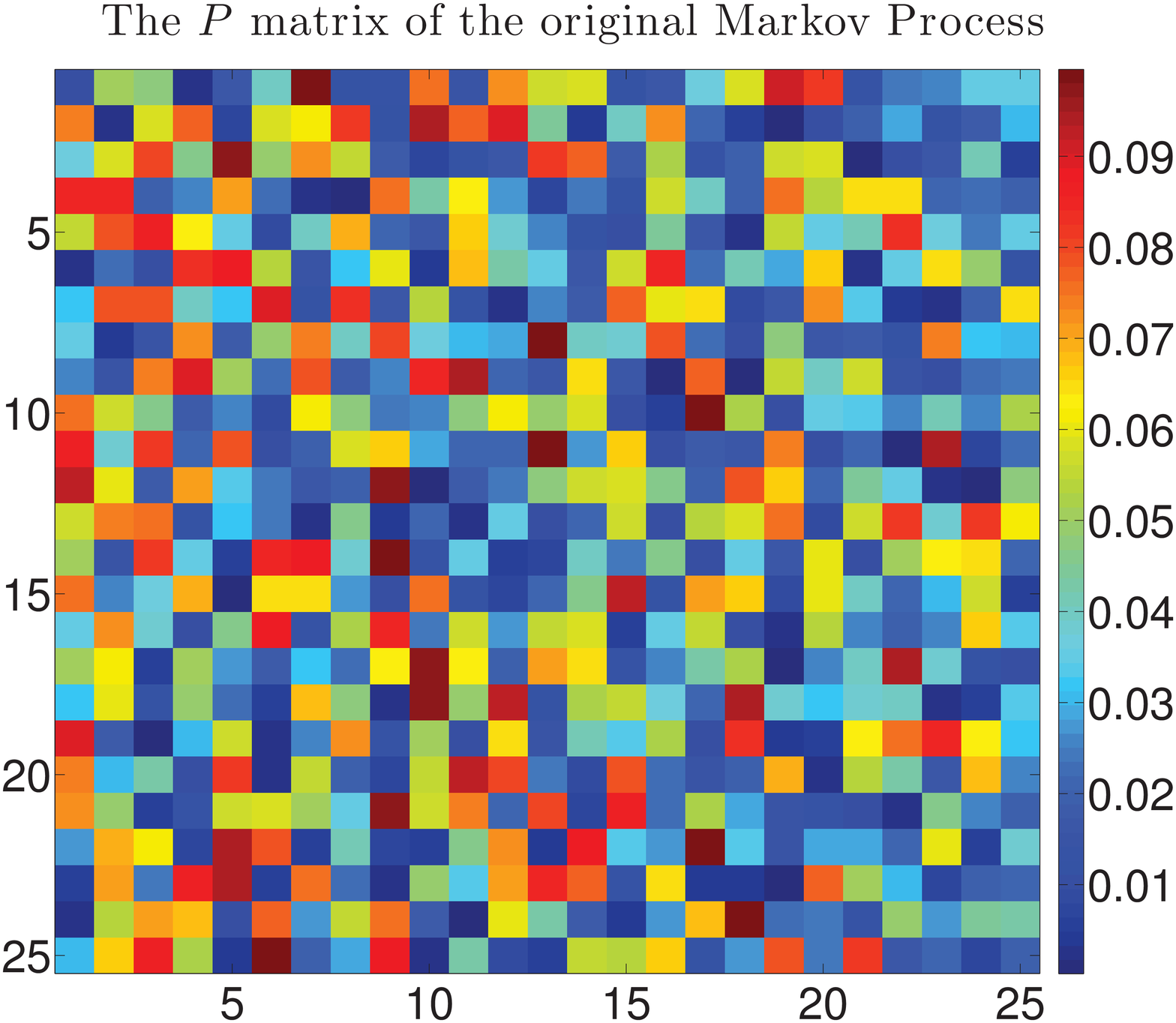}}
\subfloat[][]{
\label{fig2.2} 
\includegraphics[width=0.32\linewidth]{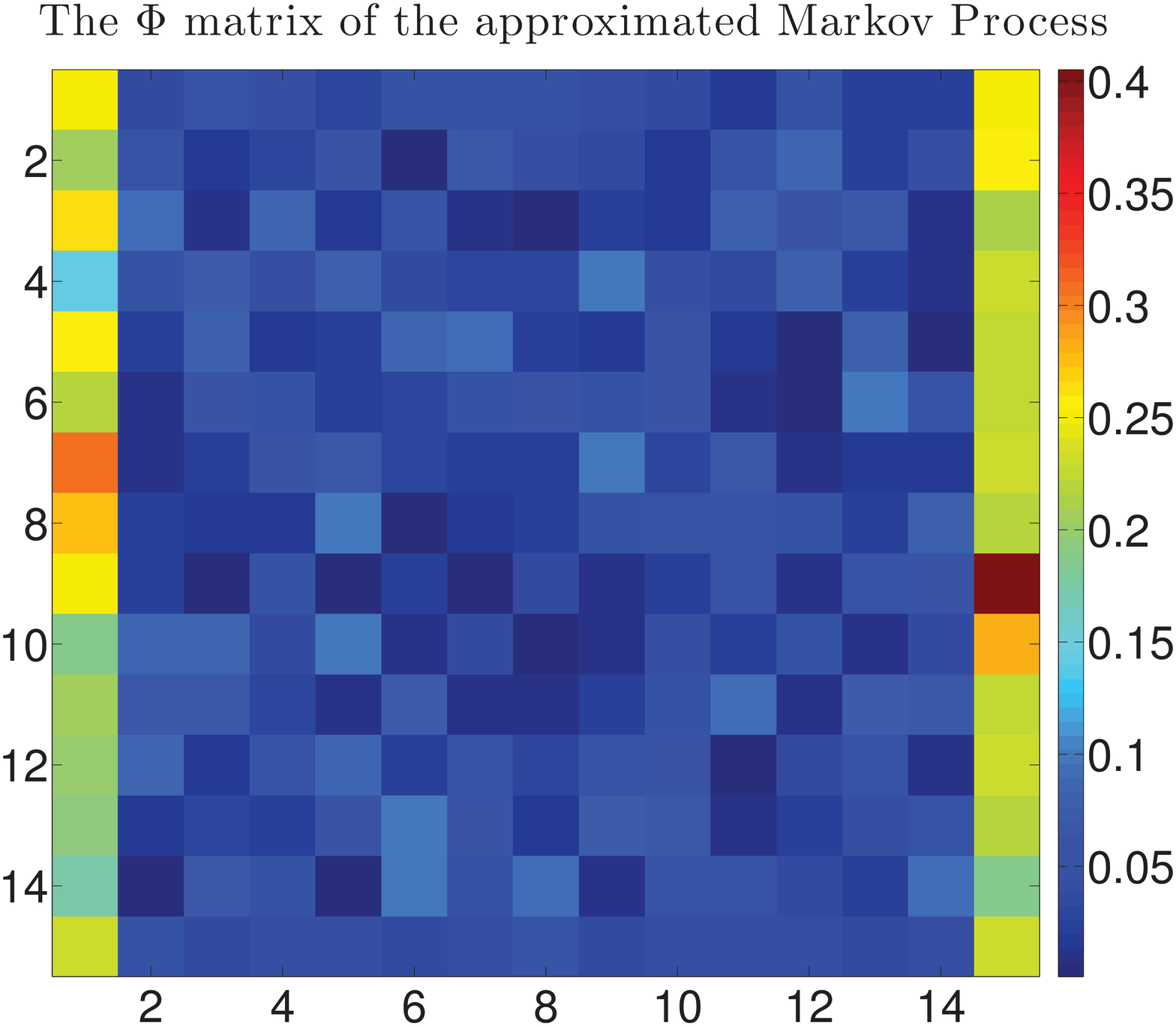}}
\subfloat[][]{
\label{fig2.3} 
\includegraphics[ width=0.32\linewidth]{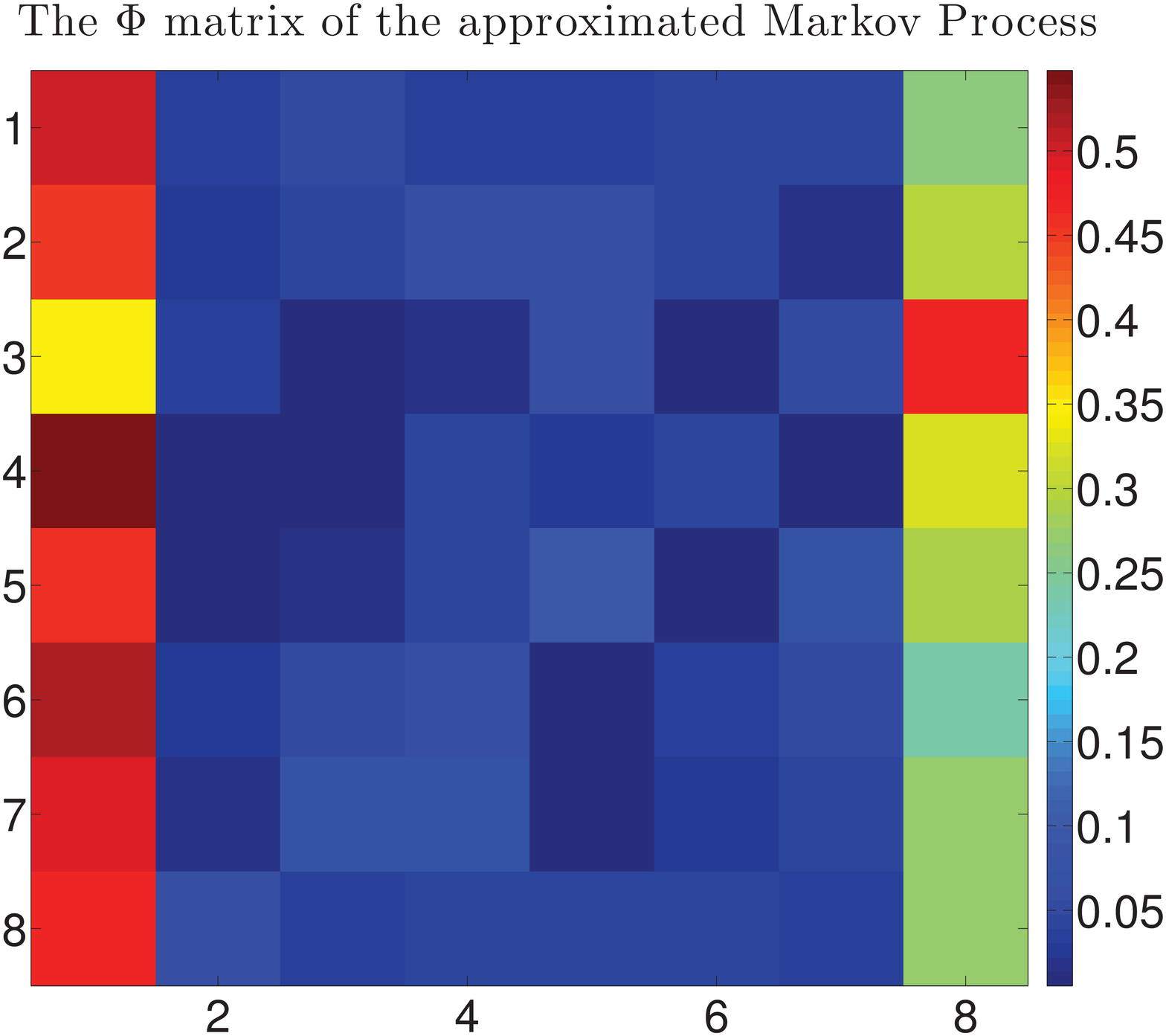}}
\hfil
\vspace{-.3cm}
\subfloat[][]{
\label{fig2.4} 
\includegraphics[ width=0.32\linewidth]{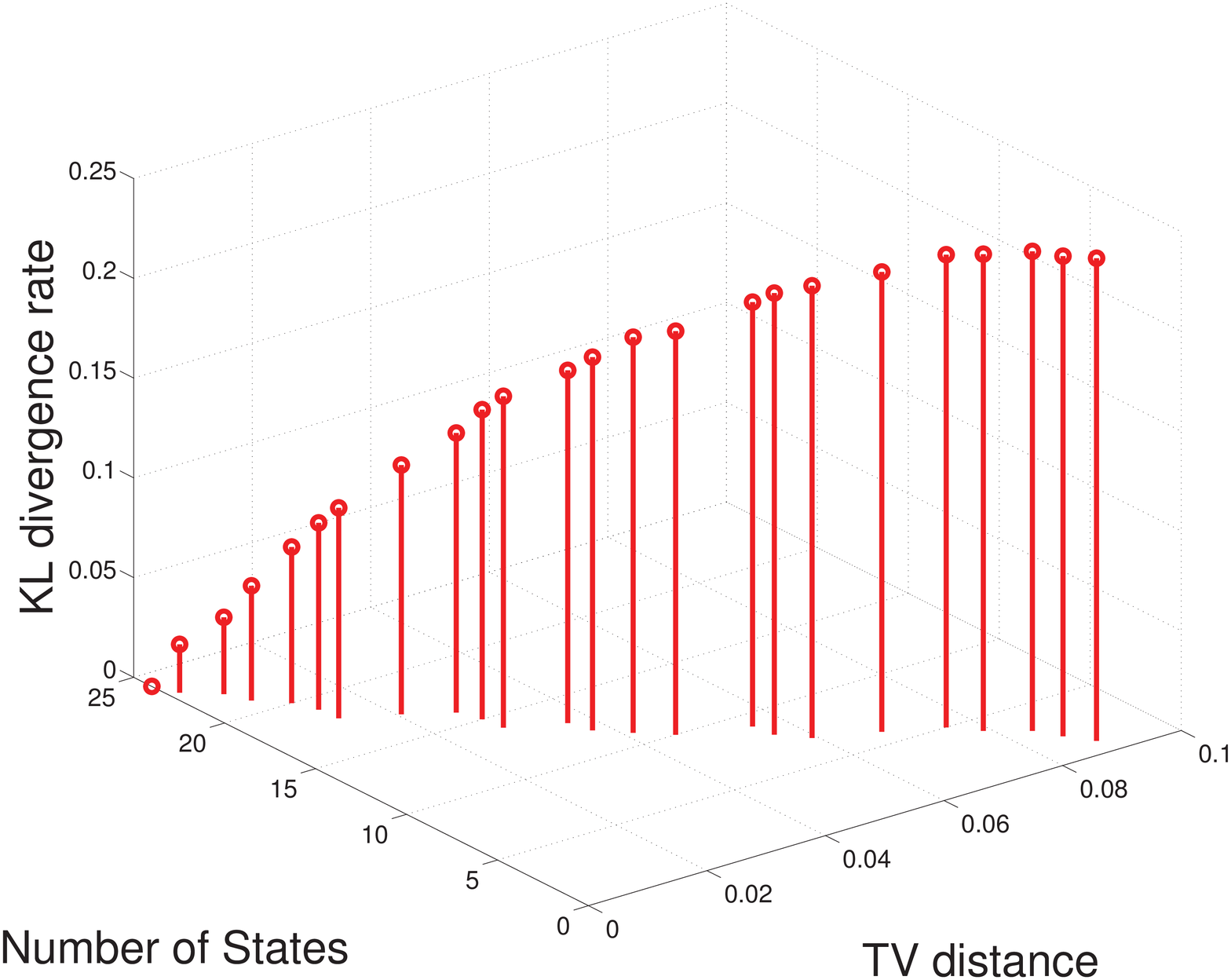}}
\subfloat[][]{
\label{fig2.5} 
\includegraphics[ width=0.32\linewidth]{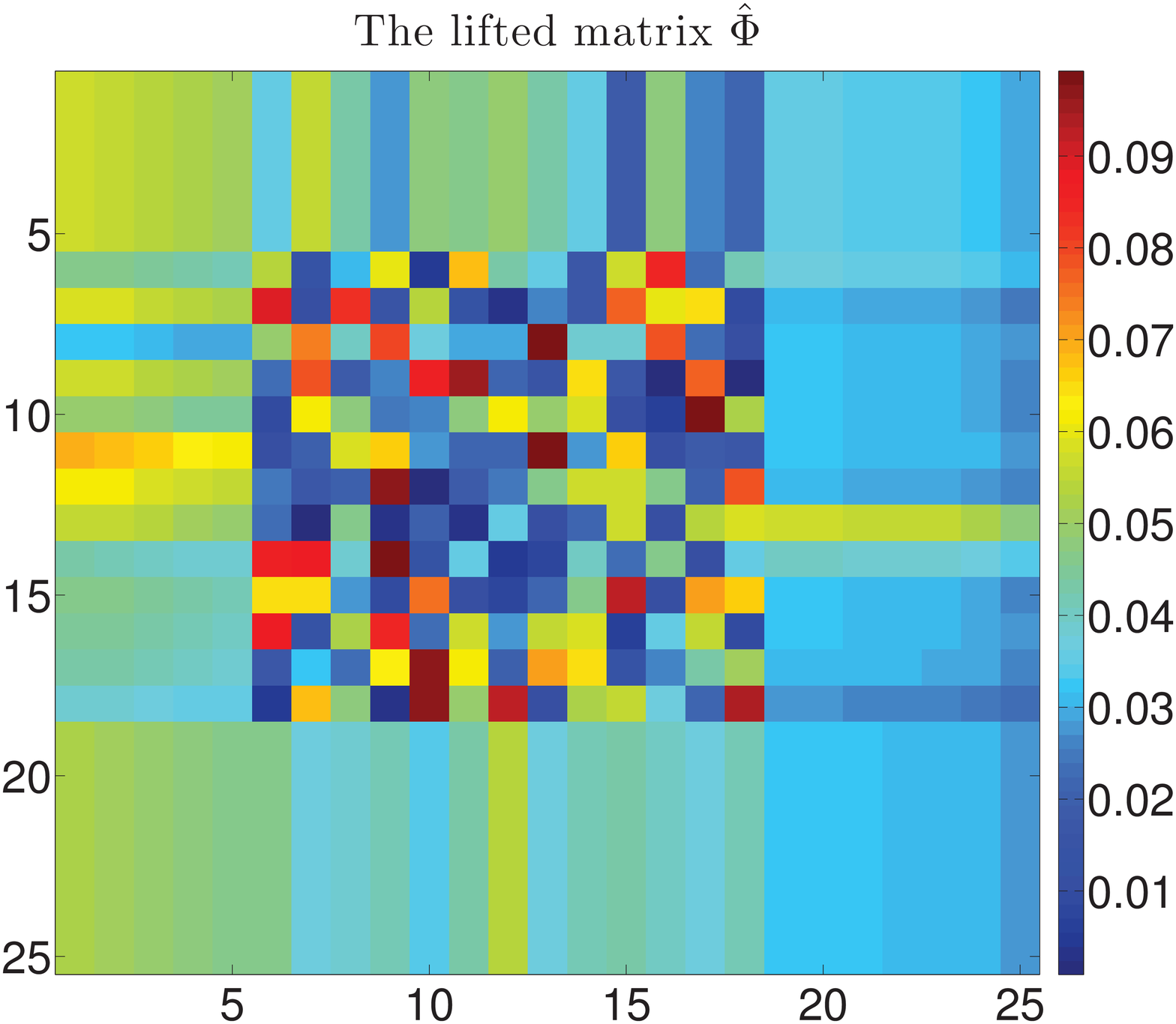}}
\subfloat[][]{
\label{fig2.6} 
\includegraphics[ width=0.32\linewidth]{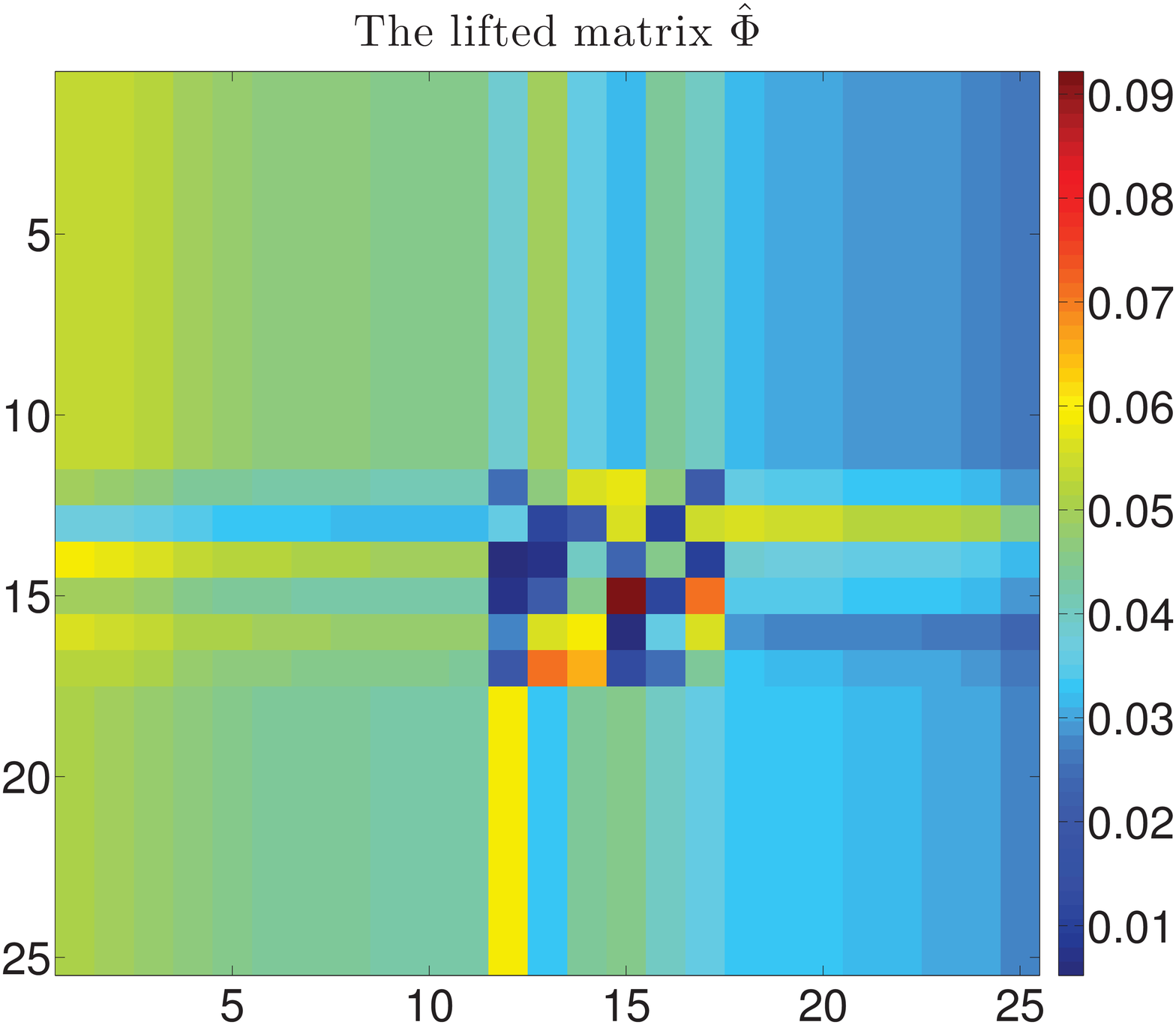}}
\caption[Optimal Solution of Example C]{Approximation results based on maximum entropy: Plot (a) depicts the $P$ matrix of the original Markov process. Plot (b) depicts a $15$-state approximation. Plot (c) depicts an $8$-state approximation. Plot (d) depicts the KL divergence rate. Plot (e) depicts the lifted $\hat{\Phi}$ matrix for the $15$-state approximation. Plot (f) depicts the lifted $\hat{\Phi}$ matrix for the $8$-state approximation.}
\label{fig2}
\end{figure*}
\subsection{Markov chain approximation based on occupancy distribution with a large number of states}\label{example2}
In this example we approximate a $25$-state Markov process based on occupancy distribution. The transition matrix $P$ of the original Markov process is as shown in Fig.\ref{fig1}\subref{fig1.1}, in which the color of the $ith$ row and $jth$ column represents the $P_{ij}$ element as indicated by the color bar. Then, based on the resulting values of $\mu_i$, $\forall i\in \calX$, the state space $\calX$ is partitioned into $16$ disjoint sets, where
\begin{IEEEeqnarray*}{c}
\calX^0{=}\{1\},~\calX_0{=}\{25\},~\calX_1{=}\{24,23\},~\calX_2{=}\{22\},~\calX_3{=}\{21\},
\calX_4{=}\{20,19\},~\calX_5{=}\{18,16\},\\ \calX_6{=}\{15\},~\calX_7{=}\{14,13\},
\calX_8{=}\{12\},~\calX_9{=}\{11,10\},~\calX_{10}{=}\{9\},~\calX_{11}{=}\{8,7\},\\ \calX_{12}{=}\{6,5\},~ \calX_{13}{=}\{4,3\},~\calX_{14}{=}\{2\}.
\end{IEEEeqnarray*}
Fig.\ref{fig1}\subref{fig1.4} depicts the KL divergence rate as a function of the number of the states of the approximated Markov process and also as a function of the TV parameter R for values where a reduction of the states occurs, due to the \textit{water-filling} behaviour of the solution. Fig.\ref{fig1}\subref{fig1.2}-\subref{fig1.5} depict the $\Phi$ matrix and the corresponding lifted matrix $\hat{\Phi}$ of the approximated Markov process, when the $25$-state Markov process is approximated by a $15$-state Markov process. Similarly, Fig.\ref{fig1}\subref{fig1.3}-\subref{fig1.6} depict $\Phi$ and $\hat{\Phi}$ when the $25$-state Markov process is approximated by an $8$-state Markov process.

\subsection{Markov chain approximation based on maximum entropy with a large number of states}\label{example3}
In this example we approximate a $25$-state Markov process based on maximum entropy. The transition matrix $P$ of the original Markov process is as shown in Fig.\ref{fig2}\subref{fig2.1}. By Remark \ref{remalg2}, the state-space $\calX$ is partitioned into $25$ disjoint sets, where $\calX^0=\{25\}$, $\calX_0=\{1\}$ and $\calX_k=\{k+1\}$ for $k=1,\dots,23$. Similarly to example \ref{example2}, Fig.\ref{fig2}\subref{fig2.4} depicts the KL divergence rate as a function of the number of the states of the approximated Markov process and as a function of TV parameter for values where an aggregation of the states occurs. It is worth noting, that the approximation based on maximum entropy principle is much faster, in terms of TV parameter, compared to the approximation based on occupancy and this is due to the \textit{water-filling like} behavior of the solution. Fig.\ref{fig2}\subref{fig2.2}-\subref{fig2.5} and \ref{fig2}\subref{fig2.3}-\subref{fig2.6} depict the $\Phi$ matrix and the corresponding lifted matrix $\hat{\Phi}$ when the original Markov process is approximated by a $15$-state and an $8$-state Markov process, respectively.

%
%
%
%
\section{Conclusion}\label{sec.conclusions}

In this work, we present two methods of approximating a FSM process by another process, with fewer states. The first method, utilizes the total variation distance to discriminate the transition probabilities of a high dimensional FSM process by a reduced order Markov process, and hence a direct method for a Markov by Markov approximation is obtained. The second method, utilizes total variation distance as a new discrepancy measure, and the problem is formulated using: (a) maximization of an average pay-off functional with respect to the approximated invariant probability, and, (b) maximization of the entropy of the approximated invariant probability, both subject to a constraint on the total variation distance metric between the invariant probability of the original Markov process and that of the approximated process. Then, by utilizing the obtained solution, we studied the problem of approximating a FSM process with another FSM process of reduced order with respect to the Kullback-Leibler divergence rate. Examples are  included to demonstrate the approximation approach for each of the two methods.

%
%
%
%
\appendix
\label{app1}
Before we proceed with the proof of Theorem \ref{MCthmocc}, we give the following Lemma in which lower and upper bounds, which are achievable, are obtained.
\begin{lemma}\label{lemmabound}
\noi 
\ben
 \item[(a)] Upper Bound. \be \sum_{j\in\calX}\ell_j\Xi_{ij}^+\mu_i\leq\ell_{\max}\left(\frac{\alpha_i\mu_i}{2}\right). \label{ub1}\ee
The bound holds with equality if 
\een
\begin{equation}\label{ub1.eq.cond}
\sum_{\mathclap{j\in\calX^0}} P_{ij}{+}\frac{\alpha_i}{2}{\leq} 1,\quad \sum_{\mathclap{j\in\calX^0}}\Xi_{ij}^+{=}\frac{\alpha_i}{2}, \quad \Xi_{ij}^+{=}0,\ \forall j{\in}\calX{\setminus}\calX^0.
\end{equation}
 \ben
\item[(b)] Lower Bound.
\ben \item[] \ben  \item[Case 1)] If $\sum_{j\in \calX_0}P_{ij}-(\alpha_i/2)\geq 0$ then \een \een \een
\be \sum_{j\in\calX}\ell_j\Xi_{ij}^-\mu_i\geq\ell_{\min}\left(\frac{\alpha_i\mu_i}{2}\right) \label{lb1}.\ee
\ben \item[] \ben
\item[] The bound holds with equality if 
\een \een
\begin{equation}\label{lb1.eq.cond}
\sum_{\mathclap{j\in \calX_0}}P_{ij}{-}\frac{\alpha_i}{2}{\geq} 0,\quad \sum_{\mathclap{j\in \calX_0}}\Xi_{ij}^-{=}  \frac{\alpha_i}{2},\quad \Xi_{ij}^-{=}0,\ \forall j{\in} \calX{\setminus}\calX_0.
\end{equation}
\ben \item[] \ben \item[] \ben \item[Case 2)] If $\sum_{s=1}^k\sum_{j\in\calX_{s-1}}P_{ij}-(\alpha_i/2)\leq 0$ for any $k\in\{1,2,\hdots,r\}$ then \een \een \een
 \begin{equation}\label{ulb1}
 \sum_{j\in\calX}\ell_j\Xi_{ij}^-\mu_i\geq \ell(\calX_k)\Big(\frac{\alpha_i\mu_i}{2}-\sum_{s=1}^k\sum_{j\in\calX_{s-1}}P_{ij}\mu_i\Big)\\+\sum_{s=1}^k\sum_{j\in\calX_{s-1}}\ell_jP_{ij}\mu_i.
 \end{equation}
\ben \item[] \ben \item[] \ben \item[] Moreover, equality holds if \een \een \een
\begin{subequations}
\label{all2}
\begin{align}
&\sum_{j\in\calX_{s-1}}\Xi_{ij}^-=\sum_{j\in\calX_{s-1}}P_{ij}, \hso\mbox{for all}\hso s=1,2,\hdots,k,\label{all2a}\\
&\sum_{j\in \calX_k}\Xi_{ij}^-=\Big(\frac{\alpha_i}{2}-\sum_{s=1}^k\sum_{j\in\calX_{s-1}}P_{ij}\Big),\label{all2b}\\
&\sum_{s=0}^k\sum_{j\in\calX_{s}}P_{ij}-\frac{\alpha_i}{2}\geq 0,\label{all2c}\\
&\Xi_{ij}^-=0 \hso\mbox{for all}\hso j\in\calX\setminus\calX_0\cup\calX_1\cup\hdots\cup\calX_k.
\end{align}
\end{subequations}
\end{lemma}

\begin{IEEEproof}
Part (a): First, we show that inequality \eqref{ub1} holds.
\begin{equation*}
\sum_{j\in\calX}\ell_j\Xi_{ij}^+\mu_i\leq \ell_{\max}\mu_i\sum_{j\in\calX}\Xi_{ij}^+=\ell_{\max}\left(\frac{\alpha_i\mu_i}{2}\right).
\end{equation*}
Next, we show that under the stated conditions \eqref{ub1.eq.cond} equality holds.
\begin{equation*}
\sum_{j\in\calX}\ell_j\Xi_{ij}^+\mu_i=\sum_{j\in\calX^0}\ell_j\Xi_{ij}^+\mu_i+\sum_{j\in\calX\setminus\calX^0}\ell_j\Xi_{ij}^+\mu_i\\
=\ell_{\max}\mu_i\sum_{j\in\calX^0}\Xi_{ij}^++\sum_{j\in\calX\setminus\calX^0}\ell_j\Xi_{ij}^+\mu_i=\ell_{\max}\left(\frac{\alpha_i\mu_i}{2}\right).
\end{equation*}
Part (b), case 1: First, we show that inequality \eqref{lb1} holds.
\begin{equation*}
\sum_{j\in\calX}\ell_j\Xi_{ij}^-\mu_i\geq \ell_{\min}\mu_i\sum_{j\in\calX}\Xi_{ij}^-=\ell_{\min}\left(\frac{\alpha_i\mu_i}{2}\right).
\end{equation*}
Next, we show that under the stated conditions \eqref{lb1.eq.cond} equality holds.
\begin{equation*}
\sum_{j\in\calX}\ell_j\Xi_{ij}^-\mu_i=\sum_{j\in\calX_0}\ell_j\Xi_{ij}^-\mu_i+\sum_{j\in\calX\setminus\calX_0}\ell_j\Xi_{ij}^-\mu_i
=\ell_{\min}\mu_i\sum_{j\in\calX_0}\Xi_{ij}^-+\sum_{j\in\calX\setminus\calX_0}\ell_j\Xi_{ij}^+\mu_i=\ell_{\min}\left(\frac{\alpha_i\mu_i}{2}\right).
\end{equation*}
Part (b), case 2: First, we show that inequality \eqref{ulb1} holds. Consider any $k\in\{1,2,\dots,r\}$. From Part (b), case 1, we have that
\begin{align*}
\sum_{j\in\calX \setminus \cup_{s=1}^k \calX_{s-1}}\ell_j\Xi_{ij}^-\mu_i&\geq  \min_{j\in \calX\setminus \cup_{s=1}^k\calX_{s-1}}\ell_j\sum_{j\in\calX\setminus\cup_{s=1}^k\calX_{s-1}}\Xi_{ij}^-\mu_i\\
&=\ell(\calX_k)\sum_{\mathclap{j\in\calX\setminus \cup_{s=1}^k\calX_{s-1}}}\Xi_{ij}^-\mu_i
=\ell(\calX_k)\Big(\sum_{j\in\calX}\Xi_{ij}^-\mu_i-\sum_{s=1}^k\sum_{j\in\calX_{s-1}}\Xi_{ij}^-\mu_i\Big).
\end{align*}
Hence,
\begin{equation*}
\sum_{j\in\calX}\ell_j\Xi_{ij}^-\mu_i-\sum_{s=1}^k\sum_{j\in\calX_{s-1}}\ell_j\Xi_{ij}^-\mu_i\\ \geq \ell(\calX_k)\Big(\frac{\alpha_i\mu_i}{2}-\sum_{s=1}^k\sum_{j\in\calX_{s-1}}P_{ij}\mu_i\Big),
\end{equation*}
which implies 
\begin{equation*}
\sum_{j\in\calX}\ell_j\Xi_{ij}^-\mu_i \geq \ell(\calX_k)\Big(\frac{\alpha_i\mu_i}{2}-\sum_{s=1}^k\sum_{j\in\calX_{s-1}}P_{ij}\mu_i\Big)+\sum_{s=1}^k\sum_{j\in\calX_{s-1}}\ell_jP_{ij}\mu_i.
\end{equation*}
Next, we show under the stated conditions \eqref{all2} that equality holds.
\begin{align*}
\sum_{j\in\calX}\ell_j\Xi_{ij}^-\mu_i&=\sum_{s=1}^k\sum_{j\in\calX_{s-1}}\ell_j\Xi_{ij}^-\mu_i{+}\sum_{j\in\calX_k}\ell_j\Xi_{ij}^-\mu_i{+}\sum_{j\in\calX\setminus\cup_{s=0}^k\calX_s}\ell_j\Xi_{ij}^-\mu_i\\
&=\sum_{s=1}^k\ell(\calX_{s-1})\sum_{j\in\calX_{s-1}}\Xi_{ij}^-\mu_i+\ell(\calX_k)\sum_{j\in\calX_k}\Xi_{ij}^-\mu_i\\
&=\sum_{s=1}^k\sum_{j\in\calX_{s-1}}\ell_jP_{ij}\mu_i+\ell(\calX_k)\Big(\frac{\alpha_i\mu_i}{2}-\sum_{s=1}^k\sum_{j\in\calX_{s-1}}P_{ij}\mu_i\Big).
\end{align*}

\end{IEEEproof}

\begin{IEEEproof}[Proof of Theorem \ref{MCthmocc}]
We provide the main steps for the derivation of Theorem \ref{MCthmocc}, since the methodology followed for solving Problem \ref{problem3} is similar to the one followed in \cite{ctlthem2013}. In particular, for a fixed $i\in\calX$, the solution of Problem \ref{problem3} is given by \eqref{mp} and \eqref{all3}, with proper substitution of $\nu^*\rightarrow \Phi^\dagger$ and $\mu\rightarrow P$. 

From \eqref{eq.payoff2}, the pay-off of Problem \ref{problem3} is given by
\begin{equation}\label{eq.posnegvar1}
\sum_{i\in \calX}\sum_{j\in \calX}\ell_j P_{ij}\mu_i+\max_{\Xi_{ij}}\sum_{i\in\calX}\sum_{j\in\calX}\ell_j\Xi_{ij}\mu_i .
\end{equation}
To maximize \eqref{eq.posnegvar1} we employ the fact that $\Xi$ is a finite signed measure satisfying \eqref{eq.posnegvar}. It is obvious that for each $i\in \calX$ an upper and a lower bound must be obtained for $\sum_{j\in\calX}\ell_j\Xi_{ij}^+\mu_i$ and $\sum_{j\in\calX}\ell_j\Xi_{ij}^-\mu_i$, respectively. Before proceeding with the derivation of the optimal transition probabilities $\Phi^\dagger$ based on upper and lower bounds, we discuss first the solution behavior in terms of the TV constraint given by \eqref{tvtvconstr}, that is $\sum_{i\in\calX}\alpha_i\mu_i\leq R$.

 Let $\alpha_i$, $\forall i\in\calX$, to be given by \eqref{eq.55d} (see \cite{ctlthem2013}, Lemma 3.1 and Corollary 3.3); then, it can be verified that for $R\leq R_{\max,i}$, $\forall i\in\calX$, the TV constraint holds with equality, and also that as $R$ increases (i.e., $R_{\max,i}\leq R\leq R_{\max,i+1}$, $\forall i,i+1\in \calX$), the TV constraint holds with inequality. However, the solution of \eqref{MCproblem1} with respect to the specific $i\in\calX$ for which $R\geq R_{\max,i}$ is constant and hence the overall solution of \eqref{MCproblem1} is not affected. Finally, for values of $R\geq R_{\max,i}$, $\forall i\in \calX$ the overall solution of Problem \ref{problem3} is constant, in particular, is equal to $\ell_{\max}$. The relation of TV constraint $\sum_{i\in\calX}\alpha_i\mu_i$ with the TV parameter $R$, is depicted in Fig.\ref{fig1110}. Next we proceed with the derivation of \eqref{eq.55}.
       \begin{figure}[h!]
\centering
\includegraphics[width=.85\linewidth]{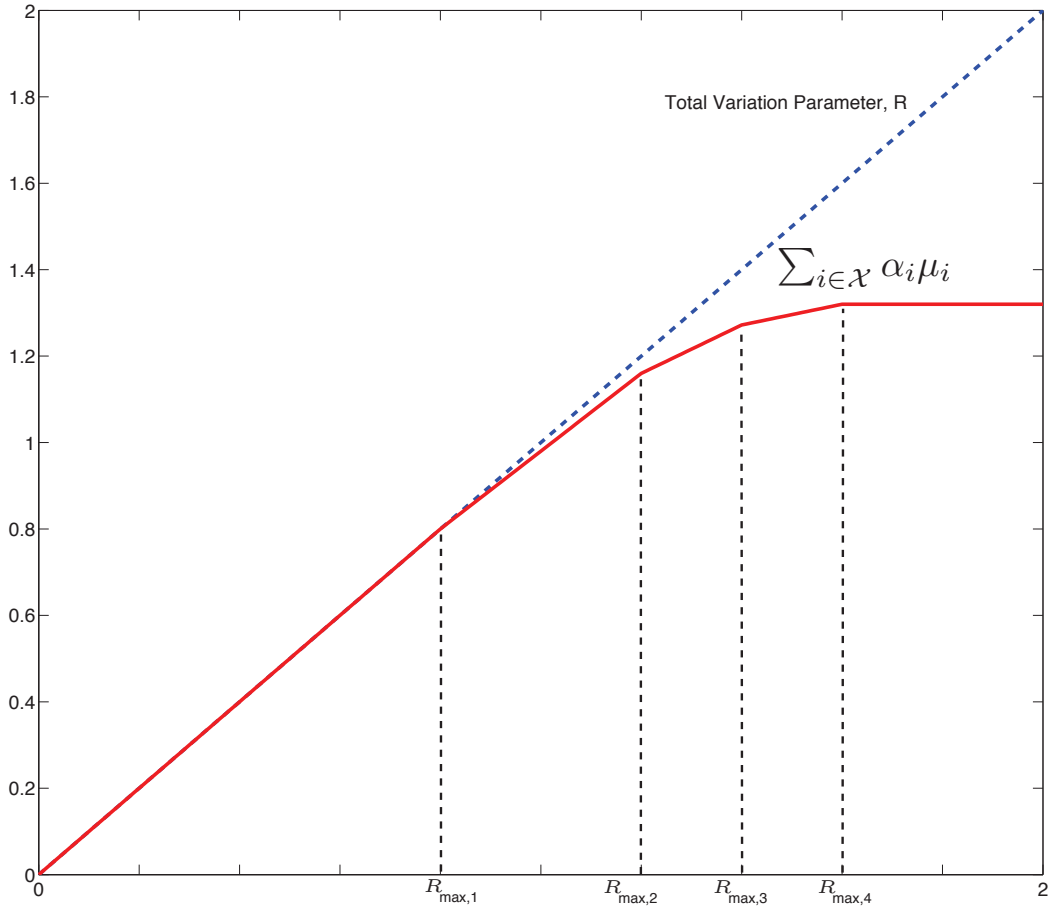}
\caption[]{Total Variation Constraint vs. Total Variation Parameter}
\label{fig1110}
\end{figure}

From Lemma \ref{lemmabound}, Part (a), the upper bound \eqref{ub1}, holds with equality if conditions given by \eqref{ub1.eq.cond} are satisfied. Note that, the first condition of \eqref{ub1.eq.cond} is always satisfied and from the second condition we have that $\sum_{j\in\calX^0}\Phi_{ij}=\sum_{j\in\calX^0}P_{ij}+\frac{\alpha_i}{2}$ and hence the optimal transition probability of each $j\in\calX^0$ is given by 
\begin{equation*}
\Phi_{ij}^\dagger=P_{ij}+\frac{\alpha_i}{2|\calX^0|},\qquad \forall j\in\calX^0.
\end{equation*}  

From  Lemma \ref{lemmabound}, Part (b), case 1, the lower bound \eqref{lb1}, holds with equality if conditions given by \eqref{lb1.eq.cond} are satisfied. Furthermore, from the second condition of \eqref{lb1.eq.cond} we have that $\sum_{j\in\calX_0}\Phi_{ij}=\sum_{j\in\calX_0}P_{ij}-\frac{\alpha_i}{2}$, and also the first condition must be satisfied, hence the optimal transition probability of each $j\in\calX_0$ is given by
\begin{equation*}
\Phi_{ij}^\dagger=\left(P_{ij}-\frac{\alpha_i}{2|\calX_0|}\right)^+,\qquad \forall j\in\calX_0.
\end{equation*}

Lemma \ref{lemmabound}, Part (b), case 1, characterize the solution for $\sum_{j\in\calX_0}P_{ij}+\frac{\alpha_i}{2}\geq 0$. Next, the characterization of solution when this condition is violated, that is, when $\sum_{s=1}^k\sum_{j\in\calX_{s-1}}P_{ij}-\frac{\alpha_i}{2}\leq 0$ for any $k\in\{1,2,\dots,r\}$ is discussed.

From Lemma \ref{lemmabound}, Part (b), case 2, the lower bound \eqref{ulb1}, holds with equality if conditions given by \eqref{all2} are satisfied. Furthermore, from \eqref{all2b} we have that 
\begin{equation*}
\sum_{j\in\calX_k}\Phi_{ij}=\sum_{j\in\calX_k}P_{ij}-\Big(\frac{\alpha_i}{2}-\sum_{s=1}^k\sum_{j\in\calX_{s-1}}P_{ij}\Big),
\end{equation*}
and conditions $\frac{\alpha_i}{2}-\sum_{s=1}^k\sum_{j\in\calX_{s-1}}P_{ij}\geq 0$ and \eqref{all2c} must be satisfied, hence the optimal transition probability of each $j\in\calX_k$ is given by 
\begin{equation*}
\Phi_{ij}^\dagger=\Big(P_{ij}-\big(\frac{\alpha_i}{2|\calX_k|}-\sum_{j=1}^k\sum_{z\in\calX_{j-1}}P_{iz}\big)^+\Big)^+.
\end{equation*}
We advice the interested reader to see \cite{ctlthem2013} for additional details concerning the steps for the solution of Problem \ref{problem3}.
\end{IEEEproof}

\bibliographystyle{IEEEtran}
\bibliography{IEEEabrv,references}

\end{document}